\documentclass[microtype]{gtpart}     
\usepackage{pinlabel}  
\usepackage{graphicx}  
\usepackage[all]{xy}

\usepackage[utf8]{inputenc}
\usepackage{psfrag}
\usepackage{hyperref}
\usepackage{tabls}

\title{Knotted Legendrian Surfaces with few Reeb chords}

\author{Georgios Dimitroglou Rizell}
\givenname{Georgios}
\surname{Dimitroglou Rizell}
\address{Department of mathematics, Uppsala University, Box 480, 751 06 Uppsala, Sweden}
\email{georgios@math.uu.se}
\urladdr{}


\arxivreference{math.SG/1102.0914}
\arxivpassword{uzxiii1}

\volumenumber{}
\issuenumber{}
\publicationyear{}
\papernumber{}
\startpage{}
\endpage{}
\doi{}
\MR{}
\Zbl{}
\received{}
\revised{}
\accepted{}
\published{}
\publishedonline{}
\proposed{}
\seconded{}
\corresponding{}
\editor{}
\version{}

\newtheorem{thm}{Theorem}[section]    
\newtheorem{lem}[thm]{Lemma}          
\theoremstyle{definition}
\newtheorem{defn}[thm]{Definition}    
\newtheorem*{rem}{Remark}             
\newtheorem{cor}[thm]{Corollary}
\newtheorem{prop}[thm]{Proposition}

\begin{document}

\begin{abstract} 
For $g>0$, we construct $g+1$ Legendrian embeddings of a surface of genus $g$ into $J^1(\R^2)=\R^5$ which lie in pairwise distinct Legendrian isotopy classes and which all have $g+1$ transverse Reeb chords ($g+1$ is the conjecturally minimal number of chords). Furthermore, for $g$ of the $g+1$ embeddings the Legendrian contact homology DGA does not admit any augmentation over $\Z_2$, and hence cannot be linearized. We also investigate these surfaces from the point of view of the theory of generating families. Finally, we consider Legendrian spheres and planes in $J^1(S^2)$ from a similar perspective.
\end{abstract}

\maketitle



\section{Introduction}

We will consider contact manifolds of the form $J^1(M)=T^*M \times \R$, where $M$ is a 2--dimensional manifold, equipped with the contact form $\alpha:=dz+\theta$. Here $\theta=-\sum_i p^i dq^i$ denotes the canonical (or Liouville) form on $T^*M$ and $z$ is the coordinate of the $\R$-factor.

An embedded surface $L \subset J^1(M)$ is called \emph{Legendrian} if $L$ is everywhere tangent to the contact distribution $\ker (\alpha)$. The \emph{Reeb vector field}, which is defined by
\[\iota_R d\alpha=0, \: \alpha(R)=1,\]
here becomes $R=\partial_z$. A \emph{Reeb chord} on $L$ is an integral curve of $R$ having positive length and both endpoints on $L$. When considering immersed Legendrian submanifolds, we say that self-intersections are zero-length Reeb chords.

We call the natural projections
\[\Pi_F\co J^1(M) \to M \times \R,\]
\[\Pi_L\co J^1(M) \to T^*M,\]
the \emph{front projection} and the \emph{Lagrangian projection}, respectively. A Legendrian submanifold $L \subset J^1(M)$ projects to an exact immersed Lagrangian submanifold $\Pi_L(L)$ in the exact symplectic manifold $(T^*M,d\theta)$. Reeb chords of $L$ correspond to self-intersections of its Lagrangian projection. 

For a generic closed Legendrian submanifold $L \subset J^1(M)$ there are only finitely many Reeb chords, each projecting to a transverse double-point of $\Pi_L(L)$ under the Lagrangian projection. We call a Legendrian satisfying this property \emph{chord generic}.

\begin{figure}[htp]
\centering
\labellist
\pinlabel $c$ at 75 38
\pinlabel $z$ at 164 73
\pinlabel $x_1$ at 142 13
\pinlabel $x_2$ at 202 35
\endlabellist
\includegraphics{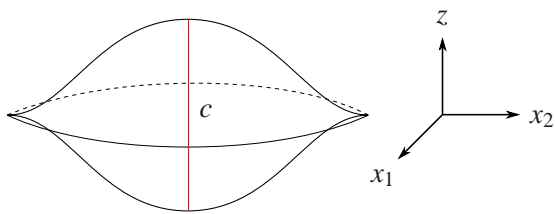}
\caption{The front projection of the standard sphere $L_{std} \subset J^1(\R^2)$.}
\label{fig:stdsphere}
\end{figure}

Let $L_{std}\subset J^1(\R^2)=\R^5$ denote the Legendrian sphere whose front projection is shown in Figure \ref{fig:stdsphere}. Note that $L_{std}$ only has one Reeb chord, and that up to isotopy it is the only known Legendrian sphere in $J^1(\R^2)$ with this property.

In Section \ref{sec:knottedsurf} we construct, for each $g>0$, the Legendrian surfaces $L_{g,k} \subset J^1(\R^2)=\R^5$ of genus $g$ by attaching $k$ ``knotted'' and $g-k$ ``standard'' Legendrian handles to $L_{std}$, where $k=0,\dots,g$. Each surface has $g+1$ transverse Reeb chords, which according to a conjecture of Arnold is the minimal number of Reeb chords for a Legendrian surface in $J^1(\R^2)$ of genus $g$. This conjecture is only known to be true for $g \leq 1$. It follows from elementary properties of generic Lagrangian immersions when $g=0$, and from Gromov's theorem of non-existence of exact Lagrangian submanifolds in $\C^n$ when $g=1$.

We will study the Legendrian contact homology of $L_{g,k}$. This theory associates a DGA (short for differential graded algebra) to a Legendrian submanifold. The DGA is then invariant up to homotopy equivalence under Legendrian isotopy. Legendrian contact homology was introduced by Eliashberg, Givental and Hofer in \cite{IntroSFT}, and by Chekanov in \cite{DiffAlg} for standard contact $\R^3$. We will also study the $L_{g,k}$ in terms of generating families (See Definition \ref{def:genfam}). We show the following theorem.

\begin{thm}
\label{thm:knottedsurfaces}
The $g+1$ Legendrian surfaces $L_{g,k} \subset J^1(\R^2)$ of genus $g$, where $k=0,\dots ,g$, are pairwise non-Legendrian isotopic. Furthermore, $L_{g,k}$ has a Legendrian contact homology DGA admitting an augmentation with coefficients in $\Z_2$ if and only if $k=0$. Also, $L_{g,k}$ admits a generating family if and only if $k=0$.
\end{thm}

\begin{rem}There is a correspondence between generating families for a Legendrian knot in $J^1(\R) \cong \R^3$ and augmentations for its DGA with coefficients in $\Z_2$. See e.g.~\cite{GeneratingFamilies}. It is not known whether a similar result holds in higher dimensions.
\end{rem}

When $k>0$, the DGA of $L_{g,k}$ with coefficients in $\Z_2$ has $1$ in the image of the boundary operator. Hence its homology vanishes, and thus it cannot be used to distinguish the different $L_{g,k}$. Moreover, it follows that its DGA has no augmentation with coefficients in $\Z_2$.

To distinguish the different $L_{g,k}$ we consider DGAs with coefficients in group ring $\Z[H_1(L_{g,k};\Z)]$ (one may also use coefficients in $\Z_2[H_1(L_{g,k};\Z)]$). We will study the augmentation varieties of these DGAs. This is a Legendrian isotopy invariant introduced by L. Ng in \cite{FramedKnot}.

In Section \ref{sec:knottedfibre} we study the following Legendrian planes. Let $F_0 := T_p^*S^2 \subset T^*S^2$ be a Lagrangian fibre and let $F_{k} \subset T^*S^2$, where $k \in \Z$, be the image of $F_0$ under $k$ iterations of a Dehn twist  along the zero-section. The plane $F_k$ coincides with $F_0$ outside of a compact set.

Since $H^1(F_k;\R)=0$, $F_k$ is an exact embedded Lagrangian submanifold and we may lift it to a Legendrian submanifold of $J^1S^2$. For the same reason, a Lagrangian isotopy of $F_k$ induces a Legendrian isotopy of the lift. Moreover, since $F_k$ is a plane, a compactly supported Lagrangian isotopy may be lifted to a compactly supported Legendrian isotopy.

By computing the Legendrian contact homology of the Legendrian lift of the link $F_{k} \cup T^*_qS^2$, we show the following.

\begin{thm}
\label{thm:knottedfibre}
There is no compactly supported Legendrian isotopy taking $F_k$ to $F_l$ if $k \neq l$. Consequently, there is no compactly supported Lagrangian isotopy taking $F_k$ to $F_l$ if $k \neq l$. However, there are such compactly supported smooth isotopies if $k \equiv l\: \mathrm{mod}2$.
\end{thm}
The effect of Dehn twists on Floer Homology was studied by P.~Seidel in \cite{ExactFloer}, and our argument is a version of it.

In Section \ref{sec:knottedsphere} we construct a Legendrian sphere $L_{knot} \subset J^1(S^2)$ with one Reeb chord which is not Legendrian isotopic to $L_{std}$. However, according to Proposition \ref{prop:ambientisotopy}, $L_{knot}$ has a Lagrangian projection which is smoothly ambient isotopic to $\Pi_L(L_{std})$. Observe that the unit disk bundle $D^*S^2$ with its canonical symplectic form is symplectomorphic to a neighbourhood of the anti-diagonal in $S^2 \times S^2$. We show the following result.

\begin{thm}
\label{thm:knottedsphere}
$L_{std}$ and $L_{knot}$ are not Legendrian isotopic. Furthermore, $\Pi_L(L_{knot})\subset D^*S^2 \subset S^2 \times S^2$ cannot be mapped to $\Pi_L(L_{std})\subset D^*S^2$ by a symplectomorphism of $S^2 \times S^2$ which is Hamiltonian isotopic to the identity.
\end{thm}

The first result is proved by computing the Legendrian contact homology of the link $L_{knot} \cup T_qS^2$. The second result follows by relating $\Pi_L(L_{knot})$ to the non-displaceable Lagrangian tori treated in \cite{ToricDegeneration}.

\section{Background}
\label{sec:background}
In this section we recall the needed results and definitions. We give a review of Legendrian contact homology, linearizations, and the augmentation variety. We also give a description of gradient flow trees, which will be used for computing the differentials of the DGAs. Finally, we briefly discuss the theory of generating families for Legendrian submanifolds.

\subsection{Legendrian contact homology}
We now recall the results in \cite{ContHomR}, \cite{ContHomP} and \cite{OrientLeg} in order to define Legendrian contact homology for Legendrian submanifolds of $J^1(M)$ with coefficients in $\Z_2$ and $\Z$. For our purposes we will only need the cases $M=\R^2$ and $M=S^2$, respectively.

The Legendrian contact homology algebra is a DGA associated to a Legendrian submanifold $L$, assumed to be chord generic, which is generated by the Reeb chords of $L$. The differential counts pseudoholomorphic disks. The homotopy type, and even the stable isomorphism type (see below), of the DGA is then invariant under Legendrian isotopy. The most obvious consequence is that the homology of the complex, the so called Legendrian contact homology, is invariant under Legendrian isotopy.

\subsubsection{The algebra}
For a chord generic Legendrian submanifold $L \subset J^1(M)$ with the set $Q$ of Reeb chords, we consider the unital algebra $\mathcal{A}_\Lambda(L)=\Lambda\langle Q \rangle$ freely generated over the ring $\Lambda$. We may always take $\Lambda=\Z_2$, but in the case when $L$ is spin we may also take $\Lambda=\Z,\Q$ or $\C$. In the latter case, the differential depends on the choice of a spin structure on $L$. For details we refer to \cite{OrientLeg}.

We will also consider the algebra $\Lambda[H_1(L;\Z)] \otimes_\Lambda \mathcal{A}_\Lambda(L)$ with coefficients in the group ring $\Lambda[H_1(L;\Z)]$.

\subsubsection{The grading}
For a Legendrian submanifold $L \subset J^1(M) = T^*M \times \R$ there is an induced \emph{Maslov class}
\[\mu\co H_1(L;\Z) \to \Z,\]
which in our setting can be computed using the following formula. Let $L$ be front generic, and let $\eta$ be a closed curve on $L$ which intersects the singular set of the front transversely at cusp edges. Recall that $z$ is the coordinate of the $\R$-factor of $J^1(M)=T^*M \times \R$. Let $D(\eta)$ and $U(\eta)$ denote the number of cusp edges transversed by $\eta$ in the downward and upward direction relative the $z$-coordinate, respectively. In \cite{NonIsoLeg} it is proved that
\begin{equation}
\label{eq:maslov}
\mu([\eta])=D(\eta)-U(\eta).
\end{equation}

We will only consider the case when the Maslov class vanishes. In this case the algebra $\mathcal{A}_\Lambda(L)$ is graded as follows. For each generator, i.e.~Reeb chord $c \in Q$, we fix a path $\gamma_c\co I \to L$ with both ends on the Reeb chord such that $\gamma_c$ starts at the point with the higher $z$-coordinate. Again, we assume that $\gamma_c$ intersects the singularities of the front projection transversely at cusp edges. We call $\gamma_c$ a \emph{capping path} for $c$. We now grade the generator $c$ by
\[|c|=\nu(\gamma_c)-1,\]
where $\nu(\gamma_c)$ denotes the \emph{Conley-Zehnder index} of $\gamma_c$.

In our setting the Conley-Zehnder index may be computed as follows. Let $f_u$ and $f_l$ be the local functions on $M$ defining the $z$-coordinates of the upper and lower sheets of $L$ near the endpoints of $c$, respectively. We define $h_c:=f_{u}-f_{l}$. Let $p\in M$ be the projection of $c$ to $M$. Observe that the fact that $c$ is a transverse Reeb chord is equivalent to $h_c$ having a non-degenerate critical point at $p$. We then have the formula
\begin{equation}
\label{eq:cz}
\nu({\gamma_c}) = D(\gamma_c)-U(\gamma_c)+\mathrm{index}_p(d^2h_c),
\end{equation}
where $D$ and $U$ are defined as above and where $\mathrm{index}_p(d^2h_c)$ is the Morse index of $h_c$ at $p \in M$. See \cite{NonIsoLeg} for a general definition of the Conley-Zehnder index and a proof of the above formula.

If the Maslov class does not vanish, we must use coefficients in $\Lambda[H_1(L;\Z)]$ to have a well-defined grading over $\Z$. Elements $A \in H_1(L;\Z)$ are then graded by
\[ |A|=-\mu(A).\]
In our cases, since $\mu$ vanishes, the coefficients have zero grading.

In the case when $L$ has several connected components, Reeb chords between two different components are called \emph{mixed}, while Reeb chords between the same component are called \emph{pure}. Mixed Reeb chords can be graded in the following way. For each pair of components $L_0, L_1$, select points $p_0 \in L_0$ and $p_1 \in L_1$ both projecting to the same point on $M$, and such that neither lies on a singularity of the front projection. Let $c$ be a mixed Reeb chord starting on $L_0$ and ending on $L_1$. A capping path is then chosen as a path on $L_1$ starting at $c$ and ending on $p_1$, together with a path on $L_0$ starting at $p_0$ and ending on $c$. The grading of a mixed chord can then be defined as before, where the Conley-Zehnder index is computed as in Formula (\ref{eq:cz}) for this (discontinuous) capping path.

Observe that the choice of points $p_0$ and $p_1$ may affect the grading of the mixed chords, hence this grading is not invariant under Legendrian isotopy in general. However, the difference in degree of two mixed chords between two fixed components is well-defined.

\subsubsection{The differential}
Choose an almost complex structure $J$ on $T^*M$ compatible with the canonical symplectic form. We are interested in finite-energy pseudoholomorphic disks in $T^*M$ having boundary on $\Pi_L(L)$ and boundary punctures asymptotic to the double points of $\Pi_L(L)$. A puncture of the disk will be called \emph{positive} in case the oriented boundary of the disk makes a jump to a sheet with higher $z$-coordinate at the Reeb chord, and will otherwise be called \emph{negative}. We assume that the chosen $J$ is \emph{regular}, i.e.~that the solution spaces of $J$-holomorphic disks with one positive puncture are transversely cut out manifolds of the expected dimension. (See \cite{ContHomP} for the existence of such almost complex structures.)

Since $\Pi_L(L)$ is an exact immersed Lagrangian, one can easily show the following formula for the (symplectic) area of a disk $D \subset T^*M$ with boundary on $\Pi_L(L)$, having the positive punctures $a_1,...,a_n$ and the negative punctures $b_1,...,b_n$:
\begin{equation}
\label{eq:area}
\mathrm{Area}(D)=\ell(a_1)+...+\ell(a_n)-\ell(b_1)-...-\ell(b_n),
\end{equation}
where $\ell(c)$ denotes the \emph{action} of a Reeb chord $c$, which is defined by
\[\ell(c):=\int_c \alpha > 0. \]
One thus immediately concludes that a non-constant pseudoholomorphic disk with boundary on $\Pi_L(L)$ must have at least one positive puncture.

Let $\mathcal{M}(a;b_1,...,b_n;A)$ denote the moduli space of pseudoholomorphic disks having boundary on $L$, a positive puncture at $a \in Q$ and negative punctures at $b_i \in Q$ in the above order relative the orientation of the boundary. We moreover require that when closing up the boundary of the disk with the capping paths at the punctures (oriented appropriately), the cycle obtained is contained in the class $A \in H_1(L;\Z)$. We define the differential on the generators by the formula
\[ \partial a = \sum_{\dim \mathcal{M}=0} |\mathcal{M}(a;b_1,...,b_n;A)|Ab_1 \cdot...\cdot b_n.\]
where $|\mathcal{M}(a;b_1,...,b_n;A)|$ is the algebraic number of elements in the compact zero-dimensional moduli space. The above count has to be performed modulo $2$ unless the moduli spaces are coherently oriented. When $L$ is spin, a coherent orientation can be given after making initial choices. If we are working with coefficients in $\Lambda$ instead of $\Lambda[H_1(L;\Z)]$, we simply project the group ring coefficient $A$ to $1$ in the above formula.

For a generic almost complex structure $J$, the dimension of the above moduli space is given by
\[ \dim \mathcal{M}(a;b_1,...,b_n;A) = |a|- |b_1| -...-|b_n|+\mu(A)-1, \]
and it follows that $\partial$ is a map of degree $-1$.

The differential defined on the generators is extended to arbitrary elements in the algebra by $\Lambda[H_1(L;\Z)]$-linearity and by the Leibniz rule 
\[ \partial(ab)=\partial(a)b+(-1)^{|a|}a\partial(b).\]

Since $L$ is an exact Lagrangian immersion, no bubbling of disks without punctures can occur, and a standard argument from Floer theory shows that $\partial^2=0$. Observe that the sum occurring in the differential always is finite because of Formula (\ref{eq:area}) and the fact that there are only finitely many Reeb chords.

\subsubsection{Invariance under Legendrian isotopy}

Let $\mathcal{A}=R\langle a_1,...,a_m\rangle$ and $\mathcal{A'}=R\langle a'_1,...,a'_m\rangle$ be free unital algebras over the ring $R$. An isomorphism $\varphi\co \mathcal{A} \to \mathcal{A}'$ of semi-free DGAs is \emph{tame} if, after some identification of the generators of $\mathcal{A}$ and $\mathcal{A}'$, it can be written as a composition of \emph{elementary} automorphisms, i.e.~automorphisms defined on the generators of $\mathcal{A}$ by
\[\varphi(a_i)=\left\{ \begin{array}{lr} a_i &  i \neq j
\\  Aa_j+\mathbf{b} & i = j \end{array}\right. \]
for some fixed $j$, where $A \in R$ is invertible, and $\mathbf{b}$ is an element of the unital subalgebra generated by $\{a_i ; \:\:\: i \neq j \}$.

The \emph{stabilization in degree $j$} of $(\mathcal{A},\partial)$, denoted by $\mathcal{S}_j(\mathcal{A},\partial)$, is the following operation. Add two generators $a$ and $b$ with $|a|=j$ and $|b|=j-1$ to the generators of $\mathcal{A}=R \langle a_1,...,a_m \rangle$. The  differential $\partial'$ of the stabilization is defined to be $\partial$ on old generators, while $\partial'(a)=b$ and $\partial'b=0$ for the new generators. It is a standard result (see \cite{DiffAlg}) that $(\mathcal{A},\partial)$ and $\mathcal{S}_j(\mathcal{A},\partial)$ are homotopy equivalent.

\begin{thm}[\cite{OrientLeg}]
Let $L \subset J^1(M)$ be a Legendrian submanifold (which is assumed to be spin and where a fixed spin-structure has been chosen in the case when $\Lambda$ has characteristic different from 2). The stable tame isomorphism class of its associated DGA $\Lambda[H_1(L;\Z)] \otimes \mathcal{A}_\Lambda(L)$ is preserved (after possibly shifting the degree of the mixed chords) under Legendrian isotopy and independent of the choice of a generic compatible almost complex structure. Hence, the homology
\[HC_\bullet(L;\Lambda[H_1(L;\Z)]):=H_\bullet (\Lambda[H_1(L;\Z)] \otimes \mathcal{A}_\Lambda(L),\partial)\]
is invariant under Legendrian isotopy. In particular, the homology
\[HC_\bullet(L;\Z_2):=H_\bullet (\mathcal{A}_{\Z_2}(L),\partial)\]
with coefficients in $\Z_2$ is invariant under Legendrian isotopy.
\end{thm}

\begin{rem}Different choices of capping paths give tame isomorphic DGAs. Changing the capping path of a Reeb chord $c \in Q$ amounts to adding a representative of some class $\eta_c \in H_1(L;\Z)$ to the old capping path. This gives the new DGA $(\Lambda [H_1(L;\Z)] \otimes_\Lambda \mathcal{A}_\Lambda,\varphi \partial \varphi^{-1})$, where
\[\varphi\co (\Lambda [H_1(L;\Z)] \otimes_\Lambda \mathcal{A}_\Lambda,\partial) \to (\Lambda [H_1(L;\Z)] \otimes_\Lambda \mathcal{A}_\Lambda,\varphi \partial \varphi^{-1}) \]
is the tame automorphism defined by mapping $c \mapsto \eta_cc$, while acting by identity on the rest of the generators.
\end{rem}
\begin{rem}
The choice of spin structure on $L$ induces the following isomorphism of the DGAs involved. Let $\mathfrak{s}_0$ and $\mathfrak{s}_1$ be two spin structures, and let $\mathfrak{s}_i$ induce the DGA $\left([H_1(L;\Z)] \otimes\mathcal{A},\partial_{\mathfrak{s}_i}\right)$. Then there is an isomorphism of DGAs (considered as $\Z$-algebras)
\[\varphi\co \left(\Z [H_1(L;\Z)] \otimes\mathcal{A},\partial_{\mathfrak{s}_0}\right) \to \left(\Z [H_1(L;\Z)] \otimes \mathcal{A},\partial_{\mathfrak{s}_1}\right)\]
defined by
\[\varphi(A) = (-1)^{d(\mathfrak{s_0},\mathfrak{s_1})(A)}A\]
for $A \in H_1(L;\Z)$, while acting by identity on all generators coming from Reeb chords. Here $d(\mathfrak{s}_0,\mathfrak{s}_1)\in H^1(L;\Z_2)$ is the difference cochain of the two spin structures.
\end{rem}

\subsubsection{Linearizations and augmentations}

Linearized contact homology was introduced in \cite{DiffAlg}. This is a stable tame isomorphism invariant of a DGA and hence a Legendrian isotopy invariant.

Let $\mathcal{A}_\Lambda=\bigoplus_{i=0}^\infty \mathcal{A}_\Lambda^i$ be the module decomposition with respect to word-length. Decompose $\partial=\oplus_i \partial_i$ accordingly and note that if $\partial_0 = 0$ on generators, it follows that $(\partial_1)^2=0$. We will call a DGA satisfying $\partial_0=0$ \emph{good}, and call the homology $H_\bullet(\mathcal{A}^1,\partial_1)$ its \emph{linearized contact homology}.

An \emph{augmentation} of $(\mathcal{A}_\Lambda,\partial)$ is a unital DGA morphism
\[\epsilon\co (\mathcal{A}_\Lambda,\partial) \to (\Lambda,0).\]
It induces a tame automorphism $\Phi^\epsilon$ defined on the generators by $c \mapsto c+\epsilon(c)$. $\Phi^\epsilon$ conjugates $\partial$ to
\[\partial^\epsilon := \Phi^\epsilon \partial (\Phi^\epsilon)^{-1}= \Phi^\epsilon \partial,\]
where $(\mathcal{A}_\Lambda,\partial^\epsilon)$ can be seen to be good. We denote the induced linearized contact homology by
\[ HLC_\bullet\left(L;\Lambda,\epsilon\right):=H_\bullet\left(\mathcal{A}(L)^1,(\partial^\epsilon)_1\right).\]

\begin{thm}[\cite{DiffAlg} 5.1]
Let $(\mathcal{A},\partial)$ be a DGA. The set of isomorphism classes of the graded vector spaces
\[ H_\bullet\left(\mathcal{A}^1,(\partial^\epsilon)_1\right)\]
for all augmentations $\epsilon$ is invariant under stable tame isomorphism. Hence, when $(\mathcal{A},\partial)$ is the DGA associated to a Legendrian submanifold, this set is a Legendrian isotopy invariant.
\end{thm}

\subsubsection{The augmentation variety}
\label{sec:augpol}

The augmentation variety was introduced in \cite{FramedKnot}. Let $\mathbb{F}$ be an algebraically closed field. In the following we suppose that $H_1(L;\Z)$ is a free $\Z$-module and that the coefficient ring $\mathbb{F}[H_1(L;\Z)]$ consists of elements of degree zero only (i.e.~that the Maslov class vanishes).

The maximal ideal spectrum
\[\mathrm{Sp}(\mathbb{F}[H_1(L;\Z)]) \simeq (\mathbb{F}^*)^{\mathrm{rank} H_1(L;\Z)}\]
can be identified with the set of unital algebra morphisms
\[\rho\co \mathbb{F}[H_1(L;\Z)] \to \mathbb{F}.\]
Extending $\rho$ by identity on the generators induces a unital DGA chain map
\[ \rho\co (\mathbb{F}[H_1(L;\Z)]\otimes \mathcal{A}_\mathbb{F},\partial) \to (\mathcal{A}_\mathbb{F},\partial^\rho := \rho\partial).\]
\begin{defn}
Let $(\mathbb{F}[H_1(L;\Z)]\otimes \mathcal{A}_\mathbb{F},\partial)$ be a DGA with coefficients in the group ring. Its \emph{augmentation variety} is the subvariety
\[\mathrm{AugVar}(\mathbb{F}[H_1(L;\Z)] \otimes \mathcal{A}_\mathbb{F},\partial)\subset \mathrm{Sp}(\mathbb{F}[H_1(L;\Z)])\]
defined as the Zariski closure of the set of points $\rho \in \mathrm{Sp}(\mathbb{F}[H_1(L;\Z)])$ for which the chain complex $(\mathcal{A}_\mathbb{F}, \partial^\rho)$ has an augmentation.
\end{defn}

This construction can be seen to be a contravariant functor from the category of finitely generated semi-free DGAs with coefficients in the group-ring $\mathbb{F} [H_1(L;\Z)]$ to the category of algebraic subvarieties of $\mathrm{Sp}(\mathbb{F}[H_1(L;\Z)])$. A unital DGA morphism will induce an inclusion of the respective subvarieties.

\begin{lem}
Let 
\[\rho\co \mathbb{F}[H_1(L,\mathbb{Z})] \to \mathbb{F}\]
be a unital algebra map, and let
\[\Phi\co(\mathbb{F}[H_1(L,\Z)] \otimes \mathcal{A}_\mathbb{F},\partial_\mathcal{A}) \to (\mathbb{F}[H_1(L,\Z)] \otimes \mathcal{B}_\mathbb{F},\partial_\mathcal{B})\]
be a unital DGA morphism. The existence of an augmentation of $\left(\mathcal{B}_\mathbb{F},\partial^\rho_\mathcal{B}\right)$ implies the existence of an augmentation of $\left(\mathcal{A}_\mathbb{F},\partial^\rho_\mathcal{A}\right)$.
\end{lem}

\begin{proof}
Augmentations pull back with unital DGA morphisms. The proposition follows from the fact that the induced map
\[\rho\Phi\co (\mathcal{A}_\mathbb{F},\partial^\rho_\mathcal{A}) \to (\mathcal{B}_\mathbb{F},\partial^\rho_\mathcal{B})\]
is a unital DGA morphism:
\[\rho \Phi \partial^\rho_\mathcal{A}=\rho \Phi \rho \partial_\mathcal{A}=\rho \Phi \partial_\mathcal{A} = \rho \partial_\mathcal{B} \Phi=\rho \partial_\mathcal{B}\rho\Phi=\partial^\rho_\mathcal{B}\rho\Phi.\]
\end{proof}

In particular, since stable tame isomorphic DGAs are chain homotopic, we have the following corollary.

\begin{cor}
\label{cor:augvar}
The isomorphism class of the augmentation variety of the DGA associated to a Legendrian submanifold is invariant under Legendrian isotopy.
\end{cor}

\subsection{Flow trees}
Our computations of the differentials of the Legendrian contact homology DGAs relies on the technique of \emph{gradient flow trees} developed in \cite{MorseFlow}. We restrict ourselves to the case $\dim M=2$.

\begin{defn}
Given a metric $g$ on $M$, a \emph{flow tree} on $L$ is a finite tree $\Gamma$ immersed by $f\co\Gamma \to M$, together with extra data, such that:
\begin{enumerate}
\renewcommand{\labelenumi}{(\alph{enumi})}
\item On the interior of an edge $e_i$, $f$ is an injective parametrization of a flow line of
\[-\nabla(h_i^\alpha -h_i^\beta),\]
where $h_i^\alpha$ and $h_i^\beta$ are two local functions on $M$, each defining the $z$-coordinate of a sheet of $L$. To the flow line corresponding to $e_i$ we associate its two \emph{1-jet lifts} $\phi^\alpha_i$, $\phi^\beta_i$, parameterized by
\[ \phi^\alpha_i(t)=\left(dh_i^\alpha(e_i(t)),h_i^\alpha(e_i(t))\right) \in L \subset J^1(M)=T^*M  \times \R,\]
\[ \phi^\beta_i(t)=\left(dh_i^\beta(e_i(t)),h_i^\beta(e_i(t))\right) \in L \subset J^1(M)=T^*M  \times \R,\]
and oriented by $-\nabla(h^\alpha_i-h^\beta_i)$ and $-\nabla(h^\beta_i-h^\alpha_i)$, respectively.
\item 
For every vertex $n$ we fix a cyclic ordering of the edges $\{e_i\}$. We denote the unique 1-jet lift of the $i$:th edge which is oriented towards (away from) the vertex $n$ by $\phi^{in,n}_i$ ($\phi^{out,n}_i$).
\item Consider the curves on $L \subset J^1(M)$ given by the oriented 1-jet lifts of the flow lines. Give the curves a cyclic order by declaring that for every vertex $n$ and edge $i$, the curve $\phi^{in,n}_i$ is succeeded by $\phi^{out,n}_{i+1}$. We require that the Lagrangian projections of the oriented 1-jet lifts in this order form a closed curve on $\Pi_L(L)\subset T^*(M)$.
\end{enumerate}
\end{defn}

If the 1-jet lifts $\phi_{i}^{in,n}$ and $\phi_{i+1}^{out,n}$ have different $z$-coordinates at the vertex $n$, we say that this is a \emph{puncture} at the vertex. The puncture is called \emph{positive} if the oriented curve jumps from a lower to a higher sheet relative the $z$-coordinate, and is otherwise called \emph{negative}.

We will also define a \emph{partial flow tree} as above, but weakening condition (c) by allowing 1-valent vertices $n$ for which $\Pi_L \circ \phi^{in,n}$ and  $\Pi_L \circ \phi^{out,n}$ differ at $n$. We call such vertices \emph{special punctures}.

As for punctured holomorphic disks with boundary on $\Pi_L(L)$, an area argument gives the following result.

\begin{lem}[Lemma 2.13 \cite{MorseFlow}]
\label{lem:pospunct}
Every (partial) flow tree has at least one positive (possibly special) puncture.
\end{lem}

\begin{lem}
\label{lem:negheight}
Suppose that a gradient flow tree has only one positive puncture. If we give each edge $e_i$ the orientation induced by its defining vector field $-\nabla(h_i^\alpha - h_i^\beta)$, where we for each edge have ordered the defining functions for the sheets such that $h_i^\alpha - h_i^\beta > 0$, then we obtain a directed tree (in particular, we claim that $h_i^\alpha-h_i^\beta \neq 0$ in the interior of $e_i$) with the following properties:
\begin{itemize}
\item Each vertex has at most one incoming edge.
\item For each vertex $n$ different from the one containing the positive puncture, if $e_{in}$ denotes the incoming edge, and $e_{out,1},...,e_{out,m}$ denote the outgoing edges, we have the inequality
\[\sum_{i=1}^m \left( h_{out,i}^\alpha(n) -h_{out,i}^\beta(n) \right) \leq h_{in}^\alpha(n)-h_{in}^\beta(n).\]
\end{itemize}
\end{lem}
\begin{proof}
To see that we get a well-defined directed tree, observe that if $h_i^\alpha-h_i^\beta = 0$ for an interior point of $e_i$, splitting the tree at this point then produces two partial flow trees, one of which has no positive punctures. This contradicts Lemma \ref{lem:pospunct}.

To prove the first of the two claims, assume that this is not the case, i.e.~that some vertex have at least two incoming edges. Split the tree somewhere at the incoming edge which is furthest away from the positive puncture. This produces two partial flow trees, one of which has no positive punctures. Again, this leads to a contradiction.

Finally, to prove the last claim, observe that if the inequality
\[\sum_{i=1}^m \left( h_{out,i}^\alpha(n) -h_{out,i}^\beta(n) \right) > h_{in}^\alpha(n)-h_{in}^\beta(n)\]
holds at a vertex $n$, then property (c) above implies that there is a positive puncture at the vertex.
\end{proof}

See definitions 3.4 and 3.5 in \cite{MorseFlow} for the notion of dimension of a gradient flow tree.

\begin{prop}[Proposition 3.14 \cite{MorseFlow}]
For a generic perturbation of $L$ and $g$, the flow trees with at most one positive puncture form a transversely cut out manifold of the expected dimension.
\end{prop}

Observe that since swallowtail singularities have codimension 2, a generic gradient flow tree will not pass through such a singularity. Moreover, the assumption $\dim M=2$ excludes more complicated singularities of the front projection.

\begin{center}
\begin{figure}[htp]
\begin{tabular}{|r|c|c|c|}
\hline
Vertex & Lagrangian projection & Front projection & Flow tree \\
\hline
$(P_1)$
&
\includegraphics{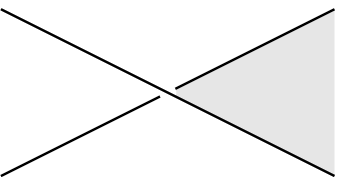}
&
\includegraphics{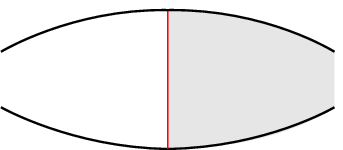}
&
\includegraphics{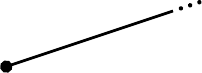}\\ \hline

$(P_2)$
&
\includegraphics{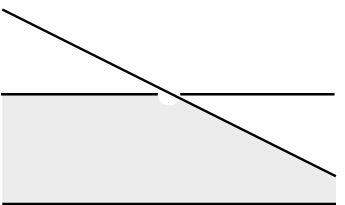}
&
\includegraphics{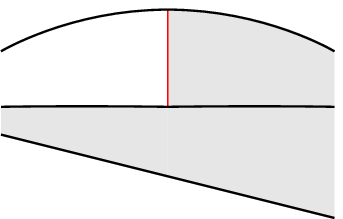}
&
\includegraphics{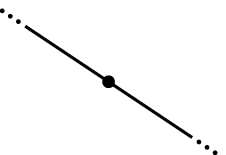}\\ \hline

$(E)$
&
\includegraphics{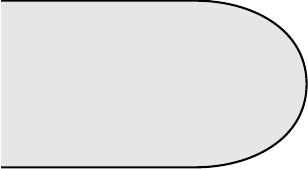}
&
\includegraphics{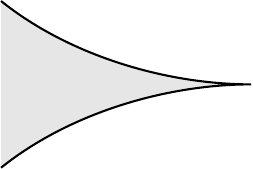}
&
\includegraphics{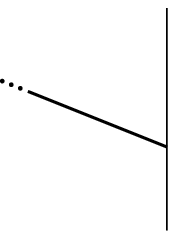}\\ \hline

$(S)$
&
\includegraphics{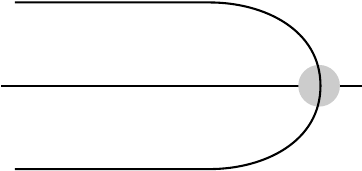}
&
\includegraphics{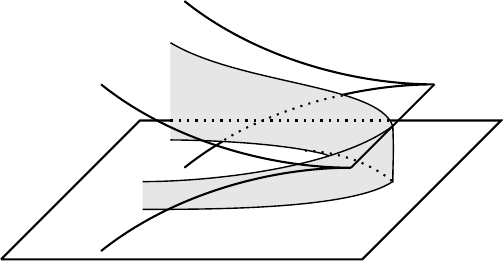}
&
\includegraphics{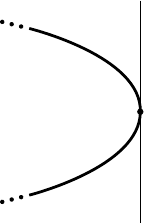}\\ \hline

$(Y_0)$
&
\includegraphics{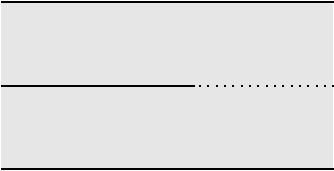}
&
\includegraphics{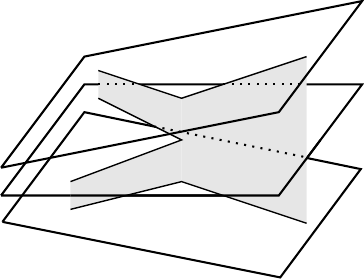}
&
\includegraphics{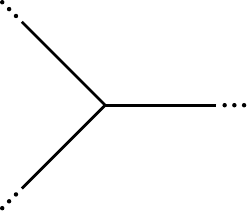}\\ \hline

$(Y_1)$
&
\includegraphics{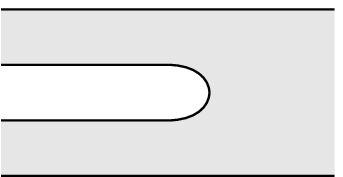}
&
\includegraphics{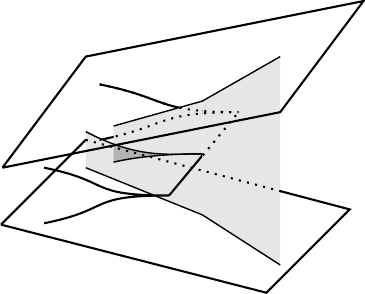}
&
\includegraphics{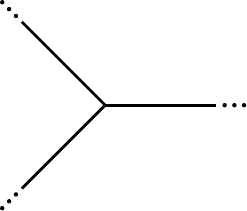}\\
\hline

\end{tabular}
\caption{$(P_1)$ and $(P_2)$ depict the punctures in the generic case. $(E)$ and $(S)$ depict the vertices corresponding to an \emph{end} and a \emph{switch}, respectively, while $(Y_0)$ and $(Y_1)$ describes the generic 3-valent vertices.}
\label{fig:trees}
\end{figure}
\end{center}

Lemma 3.7 in \cite{MorseFlow} implies that a generic, rigid, and transversely cut out gradient flow tree has no vertices of valence higher than three, and that each vertex is one of the six types depicted in Figure \ref{fig:trees}, of which the vertices $(P_1)$ and $(P_2)$ are the possible punctures. (Note that there exists both positive and negative punctures of type $(P_1)$ and $(P_2)$.)

The only picture in Figure \ref{fig:trees} which is not self explanatory is the one for $(S)$. Observe that the flow line at an $(S)$-vertex is tangent to the projection of the cusp edge. We refer to Remark 3.8 in \cite{MorseFlow} for details. 

Because of the following result, we may use rigid flow trees to compute the Legendrian Contact Homology.

\begin{thm}[Theorem 1.1 \cite{MorseFlow}]
For a generic perturbation of $L$ and the metric $g$ on $M$, there is a regular almost complex structure $J$ on $T^*M$ compatible with the canonical symplectic form, such that there is a bijective correspondence between rigid $J$-holomorphic disks with one positive puncture having boundary on a perturbation $\widetilde{L}$ of $L$, and rigid flow trees on $L$ with one positive puncture.
\end{thm}

\newpage

\subsection{Generating families}
\label{sec:generating}

\begin{defn}
\label{def:genfam}
A \emph{generating family} for a Legendrian submanifold $L \subset J^1(M)$ is a function $F\co M \times E \to \R$, where $E$ is a smooth manifold, such that
\[L=\left\{ (\mathbf{q},\mathbf{p},z) \in J^1(M);\:\:0=\frac{\partial}{\partial w^i}F(\mathbf{q},\mathbf{w}),\: p^i=\frac{\partial}{\partial q^i}F(\mathbf{q},\mathbf{w}), \: z=F(\mathbf{q},\mathbf{w}) \right\}. \]
\end{defn}

We think of $F$ as a family of functions $F_m\co E \to \R$ parameterized by $M$, and require this family to be versal. Since we are considering the case $\dim M=2$, versality implies that critical points of $F_m$ are isolated, non-degenerate outside a set of codimension 1 of $M$, possibly of $A_2$-type (birth/death type) above a codimension 1 subvariety of $M$, and possibly of $A_3$-type above isolated points of $M$.

We are interested in the case when $E$ is either a closed manifold or of the form $E=\R^N$. In the latter case we require that $F_m$ is linear and non-zero outside of a compact set for each $m \in M$.

In both cases, the Morse homology of a function $F_m$ in the family is well-defined for generic data. In the case when $E$ is a closed manifold the Morse homology is equal to $H_\bullet(E;\Z_2)$, while it vanishes in the case $E=\R^N$.

The set of generating families for $L$ is invariant under Legendrian isotopy up to stabilization of $E$ by a factor $\R^M$ and adding a non-degenerate quadratic form on $\R^M$ to $F_m$. (See \cite{CombFronts}.)

In the case $M=\R$ we have the following result. Consider the function
\[W\co M \times \R^N \times \R^N \to \R, \:\:W(m,x,y)=F(m,x)-F(m,y).\]
For sufficiently small $\delta>0$ we consider the graded vector space
\[\mathcal{G}H_\bullet(F)=H_{\bullet+N+1}(W \geq \delta,W=\delta;\Z_2).\]

In $J^1(\R)$ there are connections between generating families and augmentations. For example, we have the following result.
\begin{thm}[5.3 in \cite{GeneratingFamilies}] Let $F$ be a generic generating family for a Legendrian knot $L\subset J^1(\R)$. Then there exists an augmentation $\epsilon$ of the DGA of $L$ which satisfies
\[\mathcal{G}H_\bullet(F) \simeq HLC(L;\Z_2,\epsilon). \]
\end{thm}

\section{The front cone}
\label{sec:cone}

In this section we describe the behavior of gradient flow trees on a particular Legendrian cylinder in $J^1(\R^2)$. We need this when investigating our Legendrian surfaces, since some of them coincide with this cylinder above open subsets of $M$ in the bundle $J^1(M) \to M$. We also prove that there is no quadratic at infinity generating family for such a Legendrian surface.

\subsection{The front cone and its front generic perturbation}
We are interested in the Lagrangian cylinder embedded by
\[ S^1 \times \R \to  T^*\R^2 \simeq \C^2, \: (\theta,r) \mapsto (r+i)(\cos \theta, \sin \theta).\]
Since its Legendrian lift has a front projection given by the double cone
\[ (\theta,r) \mapsto (r\cos \theta, r\sin \theta, r) \]
in which all of $S^1\times\{0\}$ is mapped to a point, it is not front generic. We call this Legendrian cylinder the \emph{front cone}.

To make the cylinder front generic, we perturb it in the following way. Consider the plane with coordinates $\mathbf{q}=(q^1,q^2)$, and an ellipse in this plane parameterized by $\gamma(\theta)$, $\theta\in S^1$. It can be shown that
\[ F(\mathbf{q},\theta):=\|\gamma(\theta)-\mathbf{q} \| \]
is a generating family of a Legendrian cylinder for $\mathbf{q}$ in the domain bounded by the ellipse, given that its semi-axes satisfy $b\leq a<\sqrt{2}b$. When $b < a$, the front projection is generic and has four cusp edges. The projections of these correspond to points in the plane being the envelope of inward normals of the ellipse, where each normal has length equal to the curvature radius at its starting point.

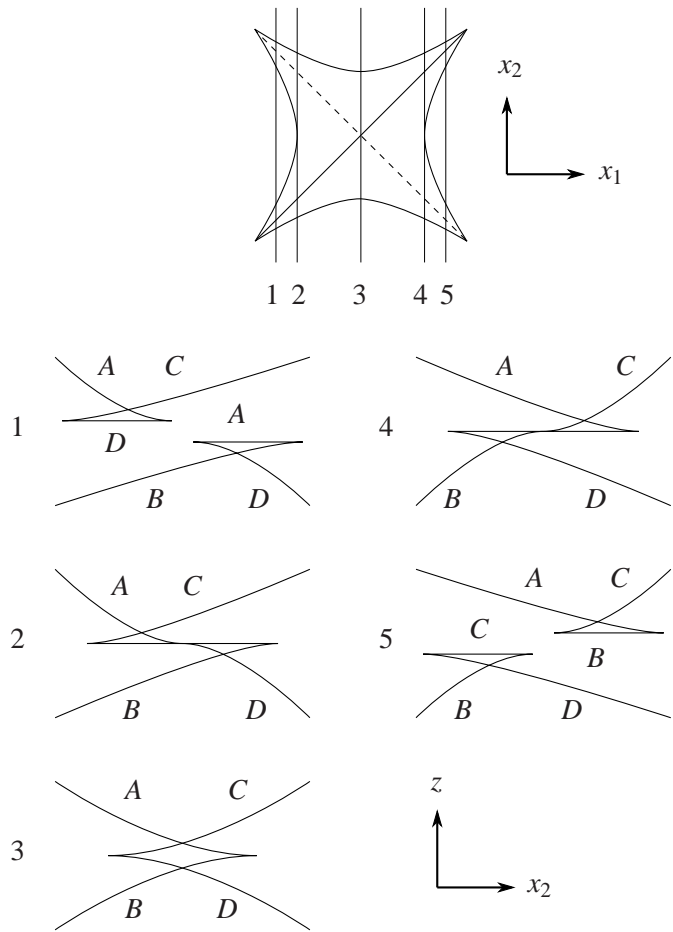
\begin{figure}[htp]
\centering
\labellist
\pinlabel $1$ at 82 240
\pinlabel $2$ at 92 240
\pinlabel $3$ at 115 240
\pinlabel $4$ at 138 240
\pinlabel $5$ at 148 240
\pinlabel $x_1$ at 210 285
\pinlabel $x_2$ at 172 325

\pinlabel $1$ at -14 190
\pinlabel $A$ at 20 213
\pinlabel $C$ at 45 213
\pinlabel $D$ at 23 184
\pinlabel $B$ at 38 163
\pinlabel $D$ at 77 163
\pinlabel $A$ at 69 195

\pinlabel $2$ at -14 110
\pinlabel $A$ at 25 130
\pinlabel $C$ at 52 130
\pinlabel $B$ at 29 83
\pinlabel $D$ at 76 83

\pinlabel $3$ at -14 30
\pinlabel $A$ at 30 53
\pinlabel $C$ at 69 53
\pinlabel $B$ at 30 8
\pinlabel $D$ at 65 8

\pinlabel $4$ at 125 190
\pinlabel $A$ at 170 213
\pinlabel $C$ at 215 213
\pinlabel $B$ at 150 163
\pinlabel $D$ at 204 163

\pinlabel $5$ at 125 110
\pinlabel $A$ at 181 133
\pinlabel $C$ at 213 133
\pinlabel $B$ at 154 83
\pinlabel $D$ at 195 83
\pinlabel $C$ at 160 114
\pinlabel $B$ at 204 103

\pinlabel $z$ at 144 55
\pinlabel $x_2$ at 183 16

\endlabellist
\includegraphics{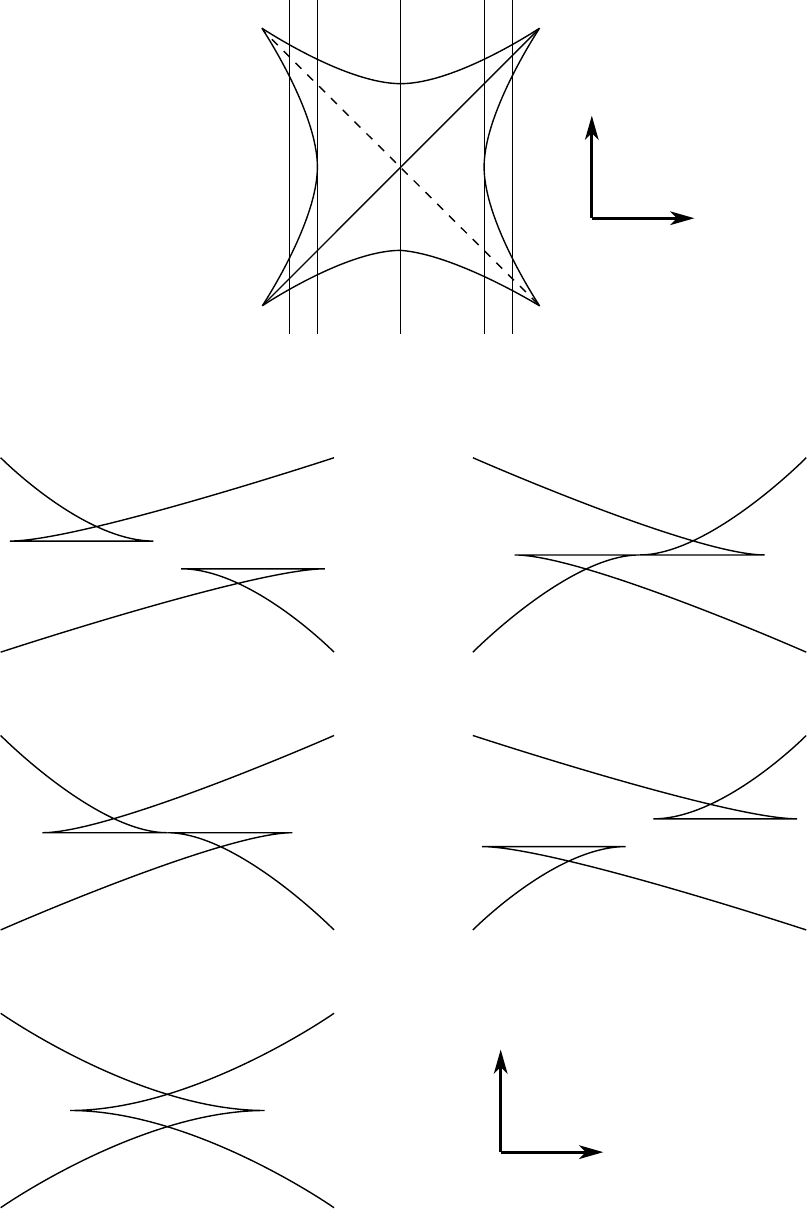}
\caption{A generic perturbation of the front cone. The top picture depicts the projection of the four cusp edges and the four swallowtail singularities; the bottom picture depicts five slices of the front projection.}
\label{fig:coneres}
\end{figure}

Degenerating the ellipse to a circle, we obtain the Legendrian cylinder corresponding to the front cone. When $b<a<\sqrt{2}b$, this ellipse is thus a generic perturbation of the front cone. Its front is described in Figure \ref{fig:coneres}.

\subsection{The gradient flow trees near the front cone}

The gradient flow outside the region in $\R^2$ bounded by the projection of the four cusp edges behaves like the gradient flow of the unperturbed cone, i.e~$-\nabla(h_u-h_l)$ points inwards to the centre, where $h_u$ is the $z$-coordinate of the upper sheet and $h_l$ is the $z$-coordinate of the lower sheet.

We will now examine the behavior of flow trees on a Legendrian submanifold $L \subset J^1(M)$ which has a front cone above some subset $U \subset M$. More precisely, we assume that after some diffeomorphism of $U \subset M$, $L$ above $U$ coincides with the front cone above an open disk centered at the origin. Moreover, we will assume that the perturbation making the front cone generic is performed in a much smaller disk.

\begin{figure}[htp]
\centering
\labellist
\pinlabel $T_1$ at 50 18
\pinlabel $T_2$ at 171 18
\endlabellist

\includegraphics{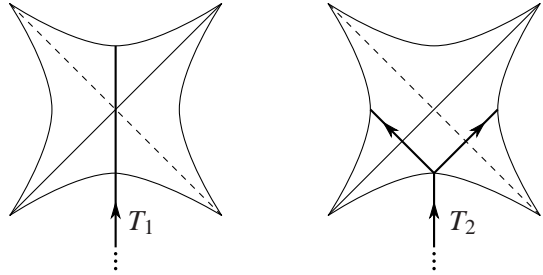}
\caption{The possibilities for a rigid gradient flow tree with one positive puncture passing through the lower cusp edge. Inside the region bounded by the cusp edges, $T_1$ is a gradient flow of the height difference of sheet $A$ and $B$. For $T_2$, the top left edge corresponds to the height difference of sheet $A$ and $D$, while the top right edge corresponds to the height difference of sheet $B$ and $C$.}
\label{fig:conedisks}
\end{figure}

\begin{prop}
\label{prop:conedisk}
Let $L \subset J^1(M)$ be a Legendrian submanifold which has a front cone above $U \subset M$, and let $T_1,\hdots, T_n$ be finitely many rigid partial gradient flow trees on $L$ with one positive puncture and which live above $M \setminus U$. Let $e_i^1,\hdots,e_i^{m_i}$ be the edges of $T_i$ which end at special punctures above $\partial U$. There is a generic resolution of the front cone with the property that if the edges $e_i^1,\hdots, e_i^{m_i}$ for each partial tree $T_i$ is to be continued to produce a rigid gradient flow tree with one positive puncture, then each edge has to be continued in one of the two ways shown in Figure \ref{fig:conedisks}.
\end{prop}

\begin{proof}
Observe that a generic resolution as shown in Figure \ref{fig:coneres} cannot have the cusp edges being tangent to $-\nabla(h^{top}-h^{bot})$ at the swallowtail singularities, where $h^{top}$ and $h^{bot}$ define the $z$-coordinates of the top and bottom sheet, respectively. A small perturbation out of the degenerate situation produces four such tangencies located arbirariy close to each of four the swallowtail singularities. These are the four points on $M$ near the cone where $(S)$-vertices might occur.

We may however perform the resolution of the front cone in such a way that when continuing each of the edges $e_i^j$ with the flow of $-\nabla(h^{top}-h^{bot})$ into the cone region, they all reach the cusp edges of the perturbed front cone at points where the difference $h^{top}-h^{bot}$ is strictly less than the absolute value of the difference in $z$-coordinate where an $(S)$-vertex may occur.

Hence, by Lemma \ref{lem:negheight}, if the completions of the $T_i$ is to have only one positive puncture it must satisfy that:
\begin{itemize}
\item Each completion of the edge $e_i^j$ cannot leave the cone region.
\item Each completion of the edge $e_i^j$ must have no $(S)$-vertex.
\end{itemize}
By dimensional reasons, if the completed gradient flow tree is to be rigid, such a completion may then neither have any $(Y_0)$-vertex in the cone region.

Thus, since the only possibility of completing the $e_i^j$ are by using $(E)$ and $(Y_1)$-vertices, we conclude that for each edge there are exactly the two possibilities analogous to the ones shown Figure \ref{fig:conedisks}.
\end{proof}

\subsection{Generating families for the front cone.}
\begin{prop}
\label{prop:conegen}
A generating family for the front cone has the property that the Morse homology of a generic function in the family is isomorphic to $H_{\bullet+i}(S^1;\Z_2)$. In particular, the front cone has no generating family with fibre $\R^N$ being linear at infinity. However, it has a generating family with fibre $S^1$.
\end{prop}
\begin{proof}
The dashed diagonal in Figure \ref{fig:coneres} corresponds to points where the $z$-coordinate of sheet $B$ is equal to that of sheet $D$. Above (below) the dashed diagonal the $z$-coordinate of sheet $B$ is greater (less) than that of sheet $D$. Similarly, the diagonal orthogonal to it corresponds to points where the $z$-coordinate of sheet $A$ is equal to that of sheet $C$. Above (below) that diagonal, the $z$-coordinate of sheet $C$ is greater (less) than that of sheet $A$.

Suppose that $F\co\R^2 \times E \to \R$ is a generating family for the front cone. We will use $A$, $B$, $C$ and $D$ to denote the critical points of the function $F_m$ in the family, where each critical point has been named after the sheet to which it corresponds. Since the pairs $AB$, $AD$, $CB$, $CD$ all cancel at birth/death singularities at the cusp edges, as seen in Figure \ref{fig:coneres}, we conclude that the critical points are graded by
\[\mathrm{index}(A)=\mathrm{index}(C)=i+1, \:\:\:\mathrm{index}(B)=\mathrm{index}(D)=i,\]
in the complex $C_\bullet(F_m;\Z_2)$, for some $i\geq0$.

For each complex, consider the pairing
\[\langle \cdot,\cdot\rangle\co C_\bullet(F_m;\Z_2) \times C_\bullet(F_m;\Z_2) \to \Z_2, \]
defined by $\langle P_i,P_j \rangle=\delta_{ij}$ for the basis $\{P_i\}$ of critical points.

Let $c \in \R^2$ denote the intersection point of the two diagonals shown in Figure \ref{fig:coneres}, i.e.~the point where the $z$-coordinates of sheets $A$ and $C$ coincide as well as those of sheets $B$ and $C$.

Suppose that $\langle\partial A,B\rangle =0$ holds for the complex $C_\bullet(F_c;\Z_2)$. Since $A$ and $B$ cancel at a birth/death singularity at the top cusp edge, there must be a handle-slide moment either from $A$ to $C$ or from $D$ to $B$ somewhere in the domain bounded by the top cusp edge and the diagonals. However, for $m$ in this domain, $F_m(A)<F_m(C)$ and $F_m(D)<F_m(B)$, so there can be no such handle slide. By contradiction, we have shown that $\langle\partial A,B\rangle=1$ must hold for $C_\bullet(F_c;\Z_2)$.

Continuing in this manner, one can show that $C_\bullet(F_c;\Z_2)$ must be the complex defined by $\partial(A)=B+D=\partial(C)$. Since the homology of this complex does not vanish, $L$ has no generating family with fibre $E = \R^N$ being linear at infinity. Assuming that $E$ is closed, we conclude that we must have $E=S^1$. Such a generating family is described in the beginning of this section.
\end{proof}

\begin{rem}
The front cone also appears twice in the conormal lift of the unknot in $\R^3$ (see \cite{InvKnot}). More precisely, let
\[f\co S^1 \hookrightarrow \R^3, \:\:f(\theta)=(\cos \theta,\sin \theta,0)\]
be the unknot. Identify $S^2$ with the unit sphere in $\R^3$. The conormal lift of $f$ is the Legendrian torus in $J^1(S^2)$ defined by the generating family
\[ F\co S^2 \times S^1 \to \R, \:\: F(\mathbf{q},\theta)=\langle \mathbf{q},f(\theta) \rangle. \]
It can be seen to have front cones above the north and south pole.
\end{rem}

\section{Knotted Legendrian surfaces in $J^1(\R^2)$ of genus $g>0$}
\label{sec:knottedsurf}

\subsection{The standard Legendrian handle}
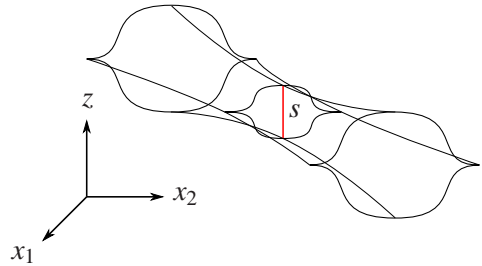
\begin{figure}[htp]
\centering
\labellist
\pinlabel $s$ at 95 49
\pinlabel $z$ at 17 53
\pinlabel $x_1$ at -7 -5
\pinlabel $x_2$ at 54 17
\endlabellist

\includegraphics{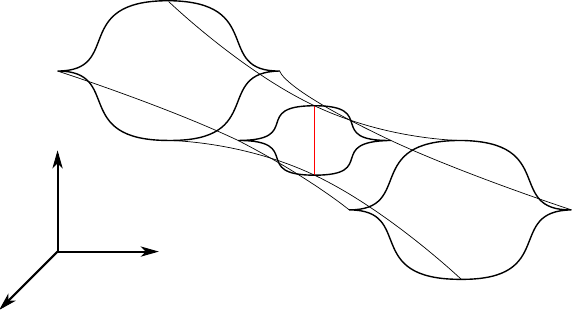}
\caption{A Legendrian handle with one saddle-type Reeb chord.}\label{fig:handle}
\end{figure}

Consider the Legendrian cylinder with one saddle-type Reeb chord whose front projection is shown in Figure \ref{fig:handle}. We will call this the \emph{standard Legendrian handle}. As described in \cite{NonIsoLeg}, we may attach the handle to a Legendrian surface by gluing its ends to cusp edges in the front projection. This amounts to performing ``Legendrian surgery'' in the case when the cusp edges are located on the same connected component, and ``Legendrian connected sum'' when the cusp edges are located on different connected components.

\begin{figure}[htp]
\centering
\labellist
\pinlabel $T_{std}$ at 7 14
\pinlabel $c_M$ at 39 17
\pinlabel $c_m$ at 166 17
\pinlabel $D_1$ at 147 37
\pinlabel $D_2$ at 166 37
\pinlabel $z$ at 228 72
\pinlabel $x_1$ at 204 10
\pinlabel $x_2$ at 267 34
\endlabellist

\includegraphics{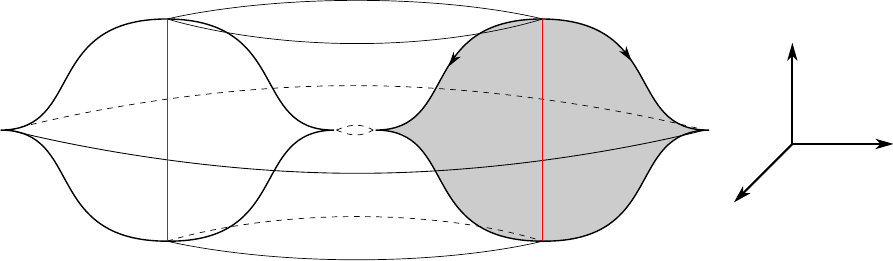}
\caption{The \emph{standard} rotational-symmetric Legendrian torus $T_{std}$ with the flow trees $D_1$ and $D_2$.}\label{fig:stdtorus}

\end{figure}
\begin{figure}[htp]
\centering

\labellist
\pinlabel $T_{knot}$ at 9 14
\pinlabel $c_M$ at 42 17
\pinlabel $c_m$ at 146 18
\pinlabel $D_1$ at 127 38
\pinlabel $D_2$ at 146 38
\pinlabel $z$ at 204 69
\pinlabel $x_1$ at 180 7
\pinlabel $x_2$ at 243 31
\endlabellist

\includegraphics{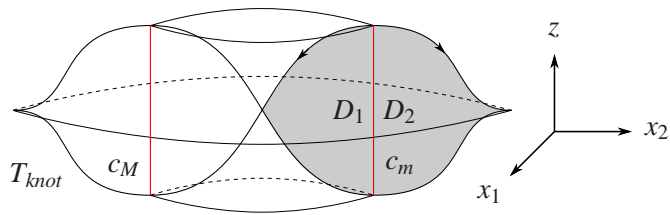}
\caption{The \emph{knotted} rotational-symmetric Legendrian torus $T_{knot}$ with the partial flow tree $D_1$ and the flow tree $D_2$.}\label{fig:knottedtorus}
\end{figure}

\subsection{Two Legendrian tori in $J^1(\R^2)$}
Consider the tori $T_{std}$ and $T_{knot}$ whose fronts are depicted in Figure \ref{fig:stdtorus} and Figure \ref{fig:knottedtorus}, respectively. For the rotational-symmetric front there is an $S^1$-family of Reeb chords. After perturbing the top sheet of the front by adding a Morse function defined on $S^1$ having exactly two critical points, we obtain two non-degenerate Reeb chords from the circle of Reeb chords: $c_M$ of maximum type and $c_m$ of saddle type. Using Formula (\ref{eq:maslov}), after making the front cone contained in $T_{knot}$ front generic as in Section \ref{sec:cone}, one sees that the Maslov class vanishes for both tori. Using the h-principle for Legendrian immersions (see \cite{Hprinciple}) one concludes that $T_{std}$ and $T_{knot}$ are regularly homotopic through Legendrian immersions.

\begin{rem}
An ambient isotopy of $\C^2$ inducing a Lagrangian regular homotopy between $\Pi_LT_{std}$ and $\Pi_LT_{knot}$ can be constructed as follows. Observe that the rotational-symmetric $\Pi_LT_{std}$ and $\Pi_LT_{knot}$ can be seen as the exact immersed Lagrangian counterparts of the Chekanov- and the Clifford torus, respectively. In the Lefschetz fibration $\C^2 \to \C$ given by $(z,w)\mapsto z^2+w^2$, with the origin as the only critical value, we get a representation of these tori as the vanishing cycle fibred over a figure-eight curve in the base. More precisely, $\Pi_L(T_{knot})$ is a figure eight curve in the base encircling the origin, while $\Pi_L(T_{std})$ is a figure eight curve not encircling the origin. From this picture it is easily seen that $\Pi_L(T_{knot})$ and $\Pi_L(T_{std})$ are ambient isotopic through immersed Lagrangians. Namely, we may disjoin the circle in the fibre from the vanishing cycle, and then pass the curve in the base through the critical value.
\end{rem}

\subsection{The Legendrian genus-$g$ surfaces $L_{g,k}$, with $0 \leq k \leq g$}
We construct the Legendrian surface $L_{g,k}$ of genus $g$ by taking the Legendrian direct sum of the standard sphere $L_{std}$ and $g-k$ copies of $T_{std}$ (or equivalently, attaching $g-k$ standard handles to $L_{std}$) and $k$ copies of $T_{knot}$. In all cases, we attach one edge of the handle to the unique cusp edge of the sphere, and the other to the outer cusp edge of the tori. The attached knotted tori will be called \emph{knotted handles}.

\begin{figure}[htp]
\centering
\labellist
\pinlabel $c_1$ at 28 58
\pinlabel $c$ at 140 58
\pinlabel $c_2$ at 271 58
\pinlabel $T_{knot}$ at 38 -5
\pinlabel $T_{std}$ at 235 -5
\pinlabel $L_{2,1}$ at 140 -9
\endlabellist
\includegraphics{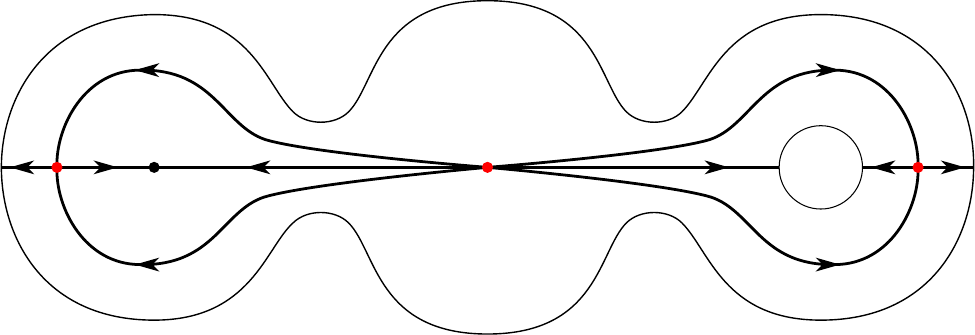}
\linebreak

\caption{The flow of $-\nabla(f_u-f_l)$, where $f_u$ and $f_l$ are $z$-coordinates of the upper and lower sheet, respectively, of the surface $L_{2,1} \simeq L_{std} \: \# \: T_{knot} \: \# \: T_{std}$.}\label{fig:surfacetop}
\end{figure}
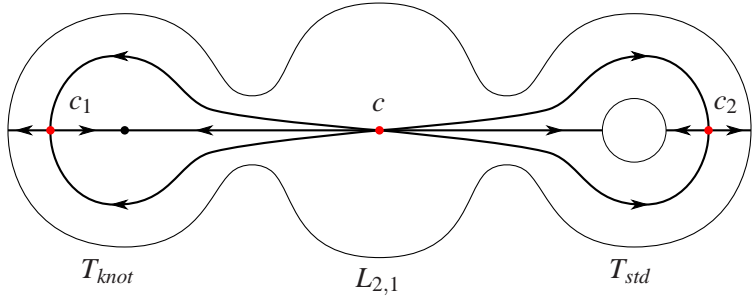

It is clear that we may cancel the Reeb chord of each standard handle connecting the sphere and the torus, which is a saddle-type Reeb chord, with the maximum-type Reeb chord $c_M$ on the corresponding torus. We have thus created a Legendrian surface having one Reeb chord $c$ coming from the sphere and Reeb chords $c_i$, for $i=1,\dots,g$, coming from the saddle-type Reeb chords on the attached tori. Such a representative is shown in Figure \ref{fig:surfacetop}.

The surfaces $L_{g,0}$ will be called the \emph{standard Legendrian surface of genus $g$} and they were studied in \cite{NonIsoLeg}. Observe that $T_{std}$ is Legendrian isotopic to $L_{1,0}$ and $T_{knot}$ is Legendrian isotopic to $L_{1,1}$.

\begin{rem}
Both $L_{1,0}$ and $L_{1,1}$ have exactly two transverse Reeb chords. This is the minimal number of transverse Reeb chords for a Legendrian torus in $J^1(\R^2)$, as follows by the adjunction formula
\[-\chi(L)+2 W(\Pi_\C(L))=0,\]
where $W$ denotes the Whitney self-intersection index. When $L$ is a torus we get that $W(\Pi_\C(L_{1,0}))=0$ which, together with the fact that there are no exact Lagrangian submanifolds in $\C^2$, implies the statement.

More generally, $L_{g,k}$ has $g+1$ transverse Reeb chords. By a conjecture of Arnold, $g+1$ is the minimal number of transverse Reeb chords for a closed Legendrian submanifold in $J^1(\R^2)$ of genus $g$.
\end{rem}

\subsection{Computation of $HC_\bullet(L_{g,k};\Z[H_1(L;\Z)])$}
We first choose the following basis $\{\mu_i,\lambda_i\}$ for $H_1(L_{g,k};\Z)$. We let $\mu_i \in H_1(L_{g,k};\Z)$ be the class which is represented by (a perturbation of) the outer cusp edge on the $i$:th torus in the decomposition
\[ L_{g,k}  \simeq L_{std} \: \#\: T_{knot} \:\#\: ... \:\#\: T_{knot}\: \#\: T_{std}\: \# \:...\: \#\: T_{std}. \]

If the $i$:th torus is standard, let $\gamma_1$ and $\gamma_2$ denote the 1-jet lifts of the flow trees on that torus corresponding to $D_1$ and $D_2$ shown in Figure \ref{fig:stdtorus}, respectively. We let $\lambda_i$ be represented by the cycle $\gamma_1-\gamma_2$.

If the $i$:th torus is knotted, $D_1$ depicted in Figure \ref{fig:knottedtorus} is a partial flow tree ending near a front cone. We choose to extend it with the partial gradient flow tree $T_1$ shown in Figure \ref{fig:conedisks} at the perturbed front cone. Again, denoting the 1-jet lifts of $D_1 \cup T_1$ and $D_2$ by $\gamma_1$ and $\gamma_2$, respectively, we let $\lambda_i$ be represented by the cycle $\gamma_1 -\gamma_2$.

\begin{lem}
$L_{g,k}$ has vanishing Maslov class, and the generators of \linebreak $\Z[H_1(L_{g,k};\Z)] \otimes \mathcal{A}(L)$ can, after an appropriate choice of capping paths, and spin stricture, be made to satisfy
\[ |c|=2, \:\:\: \partial c=0,\]
\[ |c_i|=1, \:\:\: \partial c_i=1+\lambda_i \]
if the $i$:th handle is standard, and
\[ |c_i|=1, \:\:\: \partial c_i=1+\lambda_i + \mu_i\lambda_i \]
if the $i$:th handle is knotted.
\end{lem}

\begin{proof}
Using Formula (\ref{eq:maslov}) one easily checks that the Maslov class of $L_{g,k}$ vanishes, since each closed curve on $L_{g,k}$ must traverse the cusp edges in upward and downward direction an equal number of times.

We choose the capping path for $c$ lying on $L_{std}$ (avoiding the handles) and for $c_i$ we take the 1-jet lift of $D_2$. Using Formula (\ref{eq:cz}), we compute
\[|c|=2, \:\:\: |c_i|=1. \]

After choosing a suitable spin structure on $L$ and orienting the capping operators, we get the following differential on the Reeb chords. For $c_i$ coming from a standard torus, we compute
\[ \partial c_i=1 + \lambda_i,\]
where the two terms come from the gradient flow trees $D_2$ and $D_1$, respectively. For a Reeb chord $c_i$ coming from a knotted torus we get
\[ \partial c_i=1 + \lambda_i + \mu_i \lambda_i,\]
where the first term comes from the gradient flow tree $D_2$, and the last two terms come from the partial flow tree $D_1$ approaching the front cone of the $i$:th torus. By Proposition \ref{prop:conedisk}, the edge $D_1$ can be completed in exactly two ways to become a rigid flow tree with one positive puncture: by adding the partial flow tree $T_1$ (giving the term $\lambda_i$) and by adding the partial flow tree $T_2$ (giving the term $\mu_i\lambda_i$).

For the Reeb chord $c$ coming from the sphere we get, because of the degrees of the generators, that $\partial c$ is a linear combination of the $c_i$. The relation $\partial^2=0$, together with the fact that the $\partial c_i$ form a linearly independent set, implies that
\[ \partial c=0.\]
\end{proof}

\begin{rem}
By changing the spin structure on $L_{g,k}$, we may give each term different from $1$ in the differential of a given generator an arbitrary sign.
\end{rem}

\begin{rem}For the DGA with $\Z_2$-coefficients, observe that $1$ is in the image of $\partial$ for $L_{g,k}$ when $k>0$. Hence $HC_\bullet(L_{g,k};\Z_2)=0$ in these cases.
\end{rem}

We consider the Legendrian isotopy invariant given by the augmentation variety defined in Section \ref{sec:augpol}. (In this case, it contains exactly the information given by $HC_0(L_{g,k};\Z[H_1(L;\Z)])$.)

\begin{prop}
The augmentation variety for $L_{g,k}$ over $\mathbb{C}$ is isomorphic to
\[(\C \setminus \{ 0 \})^{g-k} \times (\C\setminus \{1,0\})^k.\]
\end{prop}

\begin{proof}
After making the identification
\[\mathbb{C}[H_1(L_{g,k};\Z)] \simeq \mathbb{C}\left[Z_{\lambda_i},Z_{\lambda_i}^{-1},Z_{\mu_i},Z_{\mu_i}^{-1}\right],\]
the augmentation variety becomes
\[ \left\{1+Z_{\lambda_i}=0, \:\:\: k<i\leq g, \:\:\: 1+Z_{\lambda_j} + Z_{\lambda_j}Z_{\mu_j}= 0, \:\:\:1 \leq i \leq k \right\} \subset (\C^*)^{2g},\]
where the handles have been ordered such that the handle corresponding to $\{\lambda_i,\mu_i\}$ is standard precisely when $i>k$. To see this, observe that a DGA having no generators of degree 0 and coefficients in a field is good if and only if the differential vanishes for elements of degree 1.
\end{proof}

\begin{proof}[Proof of Theorem \ref{thm:knottedsurfaces}]
The first part follows immediately from the above proposition, together with Corollary \ref{cor:augvar}.

Since $L_{g,k}$ contains $k$ front cones, Proposition \ref{prop:conegen} gives that there does not exist any generating family for $L_{g,k}$ with vanishing Morse homology when $k>0$. Hence, there can be no generating family for $L_{g,k}$ when $k>0$, since $L_{g,k} \subset J^1(\R^2)$ is closed and a generating family for it necessarily would have vanishing Morse homology.

It can easily be checked that every $L_{g,0}$ has a linear at infinity generating families with fibre $\R$, since both the standard Legendrian handle and the standard Legendrian sphere $L_{std}$ have such generating families.
\end{proof}

\begin{rem}
Theorem \ref{thm:knottedsurfaces} also implies that the Lagrangian projections $\Pi_LL_{g,k}$ and $\Pi_LL_{g,l}$ never are Hamiltonian isotopic when $k \neq l$, even when their actions coincide.
\end{rem}

\section{Knotted Lagrangian planes in $T^*S^2$}
\label{sec:knottedfibre}
We will consider the properly embedded Lagrangian planes $F_{2k}$ being the image of the fibre $T^*_pS^2$ under the composition of $2k$ symplectic Dehn twist along the zero-section of $T^*S^2$. 

The square of a symplectic Dehn twist along the zero-section can be described by the time-$2\pi$ map of the Hamiltonian flow induced by the Hamiltonian
\[H(\mathbf{q},\mathbf{p})=\varphi(\|\mathbf{p}\|)\|\mathbf{p}\|,\]
where we have used the round metric, and where $\varphi$ is a non-decreasing function satisfying $\varphi(t)=0$ when $|t|$ is small and $\varphi(t)=1$ when $t>0$ is large. Outside of a compact set, this flow corresponds to the Reeb flow on the contact boundary $U^*S^2=\partial(D^*S^2)$ extended to a flow on the symplectization $\R \times U^*S^2 \cong T^*S^2 \setminus 0$ independently of the $\R$-factor. Consider \cite{ExactFloer} for a treatment of the symplectic Dehn twist.

\begin{figure}[htp]
\centering
\labellist
\pinlabel $D$ at 112 133
\pinlabel $b$ at 143 133
\pinlabel $a$ at 143 110
\pinlabel $F$ at 192 100
\pinlabel $F_2$ at 124 38
\endlabellist

\includegraphics{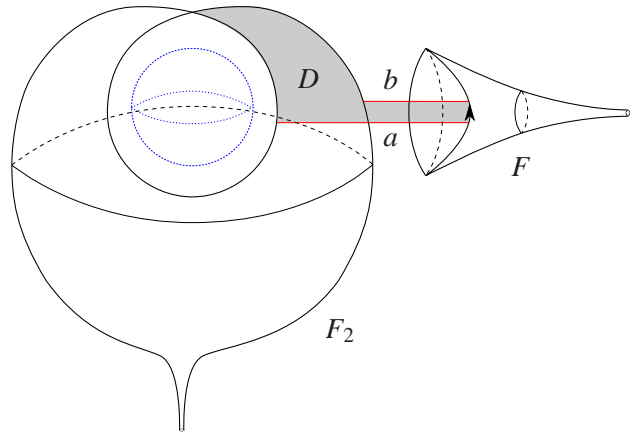}
\caption{The front projection of a Legendrian lift of $F_2 \cup F$ drawn over $S^2$, together with a partial flow tree $D$. $F$ is the Legendrian lift of a generic perturbation of $T_q^*S^2$. The innermost sphere depicts the zero-section of $T^*S^2$}\label{fig:knottedfibre}
\end{figure}

Let $q \in S^2$ be a point different from $p$. We will study the Legendrian link consisting of the Legendrian lift of $F_{2k}$ together with a Legendrian lift of a compactly supported Hamiltonian perturbation of $F:=T^*_qS^2$ as shown in Figure \ref{fig:knottedfibre}. By abuse of notation, we will sometimes use $F_{k}$ and $F$ to denote their respective Legendrian lifts in $J^1(S^2)$. We will choose the Legendrian lift of $F$ so that its $z$-coordinate is big enough to make all Reeb-chords of $F \cup F_{2k}$ start on $F_{2k}$.

\begin{figure}[htp]
\centering
\labellist
\pinlabel $F$ at 197 89
\pinlabel $F_4$ at 131 36
\pinlabel $a_1$ at 137 90
\pinlabel $a_2$ at 120 99
\pinlabel $b_1$ at 150 130
\pinlabel $b_2$ at 119 125
\endlabellist

\includegraphics{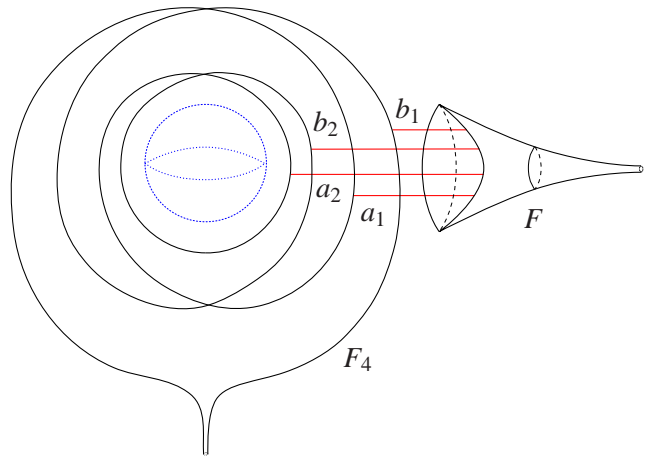}
\caption{The front projection of a Legendrian lift of $F_4 \cup F$ drawn over $S^2$.} \label{fig:knottedfibre2}
\end{figure}

The front projections of $F_2 \cup F$ and $F_4 \cup F$ are shown in Figure \ref{fig:knottedfibre} and \ref{fig:knottedfibre2}, respectively.

Since there are no pure Reeb chords, the DGA of the link is good, and since the Maslov class vanishes for both $F$ and $F_{2k}$, we may grade the DGA over $\Z$. Even though the Legendrian surfaces involved are non-compact, the Legendrian contact homology is well-defined since $F_{2k}$ and $F$ are separated by a positive distance outside of a compact set, and it is invariant under compactly supported Legendrian isotopies of $F_{2k}$. Observe that since $H^1(\R^2;\R)=0$, and since $F_{2k}$ is a plane, any compactly supported Lagrangian isotopy of $F_{2k}$ lifts to a compactly supported Legendrian isotopy.

The Legendrian link consisting of the lift of $F_{2k} \cup F$, where $F$ is translated far enough in the Reeb direction, has the mixed Reeb-chords $a_1,b_1,...,a_k,b_k$, where the Reeb chords have been labelled such that
\[ \ell(b_1) < \ell(a_1) < \ell(b_2) < \ell(a_2) ... <\ell(b_k) < \ell(a_k).\]
For $F_2 \cup F$ in Figure \ref{fig:knottedfibre}, we have $b_1=b$ and $a_1=a$.

\begin{rem}
The linearized complex for the DGA of the link is the Floer complex $CF_\bullet(F_{2k},F;\Z_2)$.
\end{rem}

\begin{lem}
\label{lem:planehom}
For a Legendrian lift of $F_{2k} \cup F$, where the lift of $F$ has been translated so that all Reeb chords start on $F_{2k}$ and end on $F$, the differential of the corresponding DGA vanishes. Consequently,
\[HLC_\bullet(F_{2k} \cup F;\Z_2) \simeq \bigoplus_{i=1}^k \Z_2 a_i \oplus \Z_2 b_i, \]
where the grading is given by
\[|b_i|=2i-1, \:\:\: |a_i|=2i,\]
and where we have chosen the unique isomorphism class of its linearized homology.
\end{lem}
\begin{proof}
We choose capping paths for each Reeb chord $c$ as follows: We fix a point $w \in F$ close to some Reeb chord endpoint. By $w' \in F_{2k}$ we denote the point on the highest sheet of $F_{2k}$ whose projection to $S^2 \subset J^1(S^2)$ coincides with that of $w$. A capping path for $c$ will be a path on $F$ starting on the endpoint of $c$ and ending at $w$, followed by a path on $F_{2k}$ starting on $w'$ and ending on the starting point of $c$. Using formula (\ref{eq:cz}) one computes
\[ |b_i|=2(i-1)+2-1, \:\:\:|a_i|=(2i-1)+2-1,\]
since the Reeb-chords all are of maximum type, and since one has to pass $2(i-1)$ (respectively $2i-1$) front cones in downward direction to go from $w'$ to $b_i$ (respectively from $w'$ to $a_i$). After perturbing each front cone as in Section \ref{sec:cone} to make it front generic, we see that traversing a front cone in downward (upward) direction amounts to traversing one cusp edge in downward (upward) direction.

By comparing indices, we immediately get that
\[ \partial_1 b_i=m_{i-1} a_{i-1}, \:\:\: \partial_1 a_i=n_i b_i, \:\:\: m_i,n_i \in \Z_2, \]
and that $\partial b_1=0$. We want to show that $m_i=n_i=0$ for all $i$. We show the case $\partial a=0$ for $F_2$ and note that the general case is analogous.

Consider the front projection of the Legendrian lift of the link shown in Figure \ref{fig:knottedfibre}. We will compute $\partial a$ by counting rigid flow trees. We are interested in flow trees having a positive puncture at $a$ and a negative puncture at $b$.

Observe that since the puncture at $b$ is negative with a maximum-type Reeb chord, it must be of type $(P_2)$. This vertex is 2-valent, with one of the edges connected to it being a flow line for the height difference of the lowest sheet of $F_2$ and a sheet of $F$, while the other edge is living on $F_2$.

Suppose that the first edge does not originate directly from a positive $(P_1)$ puncture at $a$. Thus, the edge has to end in a $(Y_0)$-vertex. However, there can be no such vertex on this edge, since this would contradict the rigidity of the flow tree.

The other edge adjacent to the $(P_2)$ vertex is a flow line living on $F_2$ approaching the front cone. Hence, we are in the situation depicted by the partial front tree $D$ in Figure \ref{fig:knottedfibre}. Applying Proposition \ref{prop:conedisk}, we conclude that this edge can be completed to a rigid flow tree with one positive puncture in exactly two ways. We have thus computed
\[\partial a=0.\]
\end{proof}

\begin{prop}
$F_{2k}$ is smoothly isotopic to $F_0$ by an isotopy having compact support.
\end{prop}
\begin{proof}
This follows from the fact that the square of a Dehn twist is smoothly isotopic to the identity by a compactly supported isotopy.
\end{proof}

\begin{proof}[Proof of Theorem \ref{thm:knottedfibre}]
Suppose that $F_{2k}$ and $F_{2l}$ are Legendrian isotopic by an isotopy having compact support, where $k,l \geq 0$. After translating $F$ far enough in the $z$-direction, we get that $F_{2k} \cup F$ is Legendrian isotopic to $F_{2l} \cup F$ by a compactly supported isotopy. Hence,
\[ HLC_\bullet(F_{2k} \cup F; \Z_2) \simeq HLC_\bullet(F_{2l} \cup F;\Z_2) \]
and we get that $k=l$ by the previous lemma. 

After applying a Dehn twist, we likewise conclude that if $F_{2k+1}$ and $F_{2l+1}$, where $k,l \geq 0$, are Legendrian isotopic by an isotopy having compact support, then $k=l$. Observe that $F_k$ and $F_l$ cannot be isotopic when $k \not\equiv l\:\mathrm{mod}2$ because of topological reasons.

Similarly, one may define $F_{k}$ for $k<0$ by applying $k$ inverses of Dehn twists. After applying $\min(|k|,|l|)$ Dehn twists to $F_k$ and $F_l$, the above result gives that $F_k$ and $F_l$ with $k,l \in \Z$ are Legendrian isotopic by a compactly supported isotopy if and only if $k=l$.

The above proposition similarly gives that $F_k$ and $F_l$ are smoothly isotopic by an isotopy having compact support if and only if $k \equiv l\:\mathrm{mod} 2$.
\end{proof}

The above computations are closely related to the result in \cite{ExactFloer}, where the existence of the following exact triangle for Floer homology is proved:
\[\xymatrix@C=-3.5em{ HF(L_0,L_1) \ar[rr] & & HF(L_0,\tau_LL_1) \ar[dl]\\
& \ar[ul] HF(L_0,L) \otimes HF(L,L_1). &
 }\]
Here $L_0, L_1$ are closed exact Lagrangians in a Liouville domain, $L$ is a Lagrangian sphere and $\tau$ is the Dehn twist along $L$ for some choice of an embedding of $L$. The map $\swarrow$ is given by the pair of pants coproduct composed with the isomorphism $HF(\tau_LL,\tau_LL_1)\simeq HF(L,L_1)$, while $\nwarrow$ is given by the pair of pants product.

\section{A knotted Legendrian sphere in $J^1(S^2)$ which is smoothly ambient isotopic to the unknot}
\label{sec:knottedsphere}
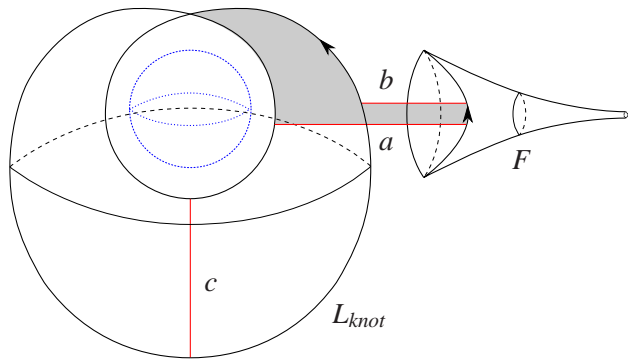
\begin{figure}[htp]
\centering
\labellist
\pinlabel $F$ at 193 75
\pinlabel $L_{knot}$ at 132 16
\pinlabel $c$ at 76 28
\pinlabel $b$ at 142 105
\pinlabel $a$ at 142 81
\endlabellist
\includegraphics{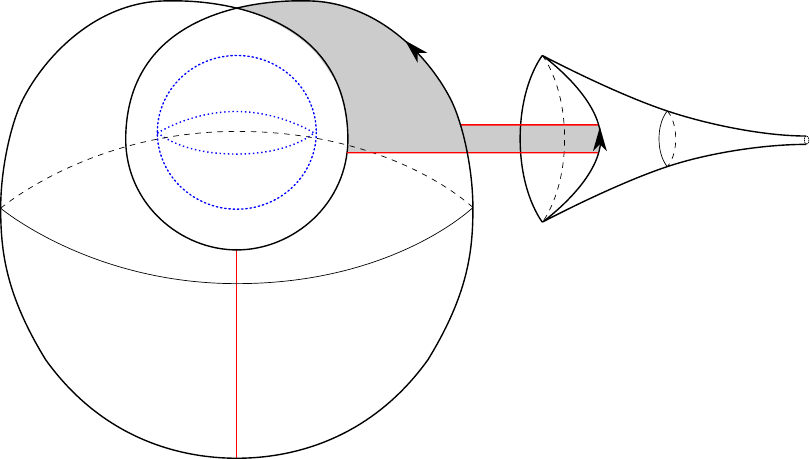}
\caption{The front projection of $L_{knot} \cup F$. The innermost sphere depicts the zero-section of $T^*S^2$.}
\label{fig:knottedsphere}
\end{figure}

Consider the front over $S^2$ given in Figure \ref{fig:knottedsphere}. It represents a Legendrian link $L_{knot} \cup F$, where $L_{knot}$ is a sphere with one maximum type Reeb chord $c$, having a rotational symmetric front, and where $F$ is the Legendrian lift of a compactly supported perturbation of a fibre $T^*_qS^2$. The Legendrian lift of the fibre has been chosen so that its $z$-coordinate is strictly larger than $\max_zL_{knot}$.

\begin{rem}
The Lagrangian projection $\Pi_L(L_{knot})$ of the sphere may be seen as the Lagrangian sphere with one transversal double point obtained by performing a Lagrangian surgery to the non-displaceable Lagrangian torus discovered in \cite{NonDisp}, exchanging a handle for a transverse double point. More precisely, the torus in $T^*S^2$ is given by the image of the geodesic flow of a fibre of $U^*S^2$. We may lift a neighbourhood of the fibre in the torus to a Legendrian submanifold in $J^1(S^2)$. The front projection of this lift looks like the front cone described in Section \ref{sec:cone}. $\Pi_L(L_{knot})$ is obtained from the torus by replacing the Lagrangian projection of the front cone with the Lagrangian projection of the two-sheeted front consisting of the graphs of functions on the form
\[f_1(x_1,x_2)=x_1^2 + x_2^2+C,\]
\[f_2(x_1,x_2)=-x_1^2-x_2^2-C,\]
where $C>0$. This produces the Reeb chord $c$ as shown at the south pole in Figure \ref{fig:knottedsphere}.
\end{rem}

\begin{rem}
One can also obtain $L_{knot}$ by the following construction, which involves a Dehn twist. Consider the standard sphere $L_{std}$ shown in Figure \ref{fig:stdsphere}, and suppose that $\Pi_L(L_{std}) \subset T^*D^2 \subset T^*\R^2$. We may symplectically embed $T^*D^2 \subset T^*S^2$ such that the zero-sections coincide. Perturb one of the sheets so that it coincides with a fibre of $T^*S^2$ in a neighbourhood of the double point. Removing a neighbourhood of this sheet and replacing it with its image under the square of a Dehn-twist along the zero-section (see Section \ref{sec:knottedfibre}), such that the Dehn twist has support in a small enough neighbourhood, yields $\Pi_LL_{knot}$.
\end{rem}

Since $L_{knot}$ has a front cone above the north pole, Proposition \ref{prop:conegen} implies that the fibre of a generating family must be $S^1$. However, the following holds.

\begin{prop}
$L_{knot}$ has no generating family
\[ F\co S^2 \times S^1 \to \R.\]
\end{prop}
\begin{proof}
Suppose there is such a generating family. Then
\[\Pi_LL_{knot}=\pi_{T^*S^2}\left(dF \cap \left( T^*S^2 \times 0_{S^1}\right)\right) \subset T^*S^2,\]
where $dF$ is considered as a section of $T^*S^2 \times T^*S^1$, $\pi_{T^*S^2}$ is the projection onto the $T^*S^2$-factor and $0_{S^1} \subset T^*S^1$ is the zero-section. Hence
\[\Pi_LL_{knot} \cap 0_{S^2} = \pi_{T^*S^2}\left(dF \cap \left(0_{S^2} \times 0_{S^1} \right)\right),\]
which by the Morse inequalities consists of at least four points when the intersection is transversal. However, one sees that $L_{knot}$ intersects the zero-section transversely in only two points, which leads to a contradiction.
\end{proof}

\begin{rem}
$L_{knot}$ can be seen to have a generating family defined on an $S^1$-bundle over $S^2$ having Euler number $1$.
\end{rem}

\begin{prop}
\label{prop:ambientisotopy}
$\Pi_L(L_{knot})$ and $\Pi_L(L_{std})$ are smoothly ambient isotopic.
\end{prop}

\begin{proof}
Consider $L_{knot}$ given by the rotation symmetric front in Figure \ref{fig:knottedsphere}. We assume that the Reeb chord $c$ is above the south pole and that the front cone is above the north pole. We endow $S^2$ with the round metric.

The goal is to produce a filling of $L_{knot}$ by an embedded $S^1$-family of disks in $T^*S^2$ with corners at $c$. More precisely, we want a map 
\[\varphi\co S^1 \times D^2 \to T^*S^2,\]
such that

\begin{itemize}
\item $\varphi$ is a diffeomorphism on the complement of $S^1 \times \{1\}$.
\item $\varphi^{-1}(c)=S^1 \times \{ 1\}$.
\item $\varphi|_{S^1 \times \partial D^2}$ is a foliation of $\Pi_L(L_{knot})$ by embedded paths starting and ending at the double point.
\item On a neighbourhood $U \supset S^1 \times \{1\}$, $\varphi|_{(\{\theta\} \times D^2) \cap U}$ maps into the plane given by $(s,t) \mapsto t\dot{\gamma}_\theta(s)\in T_{\gamma_\theta(s)}S^2 \simeq T^*_{\gamma_\theta(s)}S^2$ (identified using the round metric), where $\gamma_\theta$ is a geodesic on $S^2$ starting at the south pole with angle $\theta$.
\end{itemize}

The existence of such a filling by disks will prove the claim, since an isotopy then may be taken as a contraction of $\Pi_L(L_{knot})$ within the disks to a standard sphere contained in the neigbhourhood of $S^1 \times \{1\}$. To that end, observe that such a neighbourhood contains a standard sphere intersecting each plane $(s,t) \mapsto t\dot{\gamma}_\theta(s) \in T_{\gamma_\theta(s)}S^2$ in a figure eight curve, with the double point coinciding with that of $L_{knot}$.

We begin by considering the $S^1$-family $\{\widetilde{D}_\theta\}$ of embedded disks with boundary on $\Pi_L(L_{knot})$ and two corners at $c$, such that $\widetilde{D}_\theta$ is contained in the annulus 
\[\{t\dot{\gamma}_\theta(s) \in T_{\gamma_\theta(s)}S^2 \simeq T^*_{\gamma_\theta(s)}S^2;\:\:t \in \R, \:0\leq s<2\pi\},\]
and where $\gamma_\theta$ is the geodesic described above. (In some complex structure, these disks may be considered as pseudoholomorphic disks with two positive punctures at $c$, or alternatively, gradient flow trees on $L_{knot}$ with two positive punctures.)

Let $R_\theta\co \C^2 \to \C^2$ be the complexified rotation of $\R^2$ by angle $\theta$. We can take a chart near the north pole of $T^*S^2$ which is symplectomorphic to a neighbourhood of the origin in $\C^2 \simeq T^*\R^2$ such that:
\begin{itemize}
\item The image of $L_{knot}$ in the chart is invariant under $R_\theta$.
\item The disk in the above family corresponding to the geodesic $\gamma_\theta$ is contained in the plane $R_\theta(z,0)$.
\end{itemize}
The image of $L_{knot}$ in such a chart is shown in Figure \ref{fig:filling}.

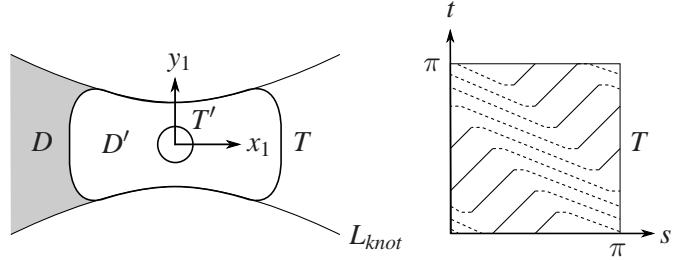
\begin{figure}[htp]
\centering

\labellist
\pinlabel $D$ at 12 36
\pinlabel $T$ at 110 36
\pinlabel $L_{knot}$ at 138 0
\pinlabel $D'$ at 40 36
\pinlabel $T'$ at 73 45
\pinlabel $x_1$ at 94 34
\pinlabel $y_1$ at 64 66

\pinlabel $T$ at 237 36
\pinlabel $s$ at 247 0
\pinlabel $t$ at 166 86
\pinlabel $\pi$ at 229 -5
\pinlabel $\pi$ at 159 65
\endlabellist

\includegraphics{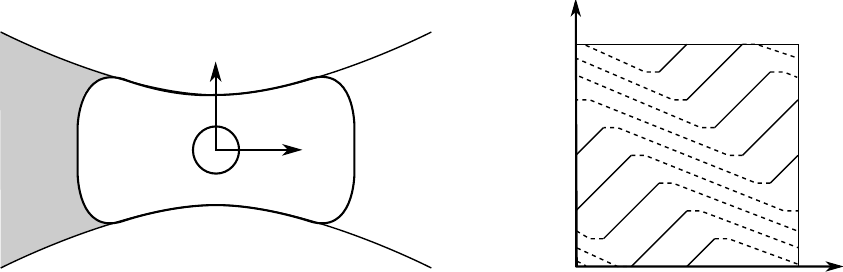}
\linebreak
\caption{A neighbourhood of the north pole of $T^*S^2$ identified as a subset of $\C^2$ where $L_{knot}$ and $T$ are invariant under $R_\theta$. The curves $(R_\theta D) \cap T$ can be completed to a foliation of $T$ by closed curves as shown on the right.}
\label{fig:filling}
\end{figure}

The torus $T$ invariant under $R_\theta$ shown in Figure \ref{fig:filling} is symplectomorphic to a Clifford torus $S^1 \times S^1 \subset \C^2$. It can be parameterized by
\[f\co \R/\pi\Z \times \R / \pi\Z \to \C^2,\]
\[(s,t) \mapsto f(s,t)=\gamma(s+t)(\cos (s-t), \sin (s-t)),\]
where $\gamma\co \R/2\pi\Z \to \C$ is a parametrization of $T \cap \left(\C \times \{0 \}\right)$ satisfying $\gamma(s+\pi)=-\gamma(s)$.

We now cut the disks in the family $\{\widetilde{D}_\theta\}$ along $T$. For each $\widetilde{D}_\theta$ we obtain three disks: $D_\theta$ coinciding with $R_\theta D$ in Figure \ref{fig:filling}, $-D_\theta$ coinciding with $R_{\theta+\pi}D$, and $R_\theta D'$.

Observe that $\{D_\theta\}$ is an $S^1$-family of embedded disks which has the right behaviour near the double-point of $\Pi_L(L_{knot})$. However, each disk in the family has a boundary arc $(R_\theta D) \cap T$ which is not on $\Pi_L(L_{knot})$. We will produce our filling $\varphi$ by gluing another family of embedded disks along these arcs.

As shown in Figure \ref{fig:filling}, each arc $(R_\theta D) \cap T$ can be extended by a curve inside $T \cap \Pi_L(L_{knot})$ to become a unique leaf in a foliation of $T$ by closed curves. This foliation extends to a filling of $T$ by an embedded $S^1$-family of disks whose interiors are disjoint from $\Pi_L(L_{knot}) \cup \bigcup_{\theta \in S^1}D_{\theta}$. Gluing these disks to the disks in $\{D_\theta\}$ will produce the required filling $\varphi$.

To see the filling of $T$ one can argue as follows. The foliation of $T$ is isotopic to a foliation where all leaves are of the form $t \equiv C$. Using this isotopy, we may create a family of annuli with one boundary component being a leaf in the foliation of $T$, and the other boundary component being the curve $\epsilon e^{i(s+C)}(\cos(s-C),\sin(s-C))$ for some $C$ and $\epsilon >0$. The latter curve is a leaf of a foliation of the torus parameterized by $\epsilon e^{i(s+t)}(\cos(s-t),\sin(s-t))$. The smaller torus, which is depicted by $T'$ in Figure \ref{fig:filling}, is clearly isotopic to $T$ by an isotopy preserving each $R_\theta D'$. Finally, the leaf in the foliation on $T'$ corresponding to $t \equiv C$ is bounded by a disk contained in the plane
\[\left(z+\frac{\epsilon}{2}e^{i2C},-iz-\frac{\epsilon}{2i}e^{i2C}\right), \: z\in \C.\]
\end{proof}

We will now compute the linearized Legendrian contact homology of the link $L_{knot} \cup F$ with coefficients in $\Z_2$. Observe that the DGA of each component is good since, as we shall see, the differential vanishes. Observe that even though $F$ is not compact, the Legendrian contact homology is still well-defined under Legendrian isotopy of the component $L_{knot}$.

\begin{lem}
The differential vanishes on all generators of the DGA for the link $L_{knot} \cup F$, and hence
\[HLC_\bullet(L_{knot} \cup F;\Z_2) \simeq \Z_2 c \oplus \Z_2 a \oplus \Z_2 b,\]
where we have chosen the unique isomorphism class of its linearized homology.
\end{lem}

\begin{proof}
Since both components have zero Maslov number, we may consider DGAs graded over $\Z$.

Formula (\ref{eq:cz}) gives $|c|=2$. We may assume that $\ell(c) < \min(\ell(a),\ell(b))$. Thus, by the area formula (\ref{eq:area}), we immediately compute $\partial(c)=0$ for the only pure Reeb chord. Hence there is an unique augmentation of the DGA for each component, namely the trivial one.

For the Reeb chords $a$ and $b$, we grade them as in the proof of Lemma \ref{lem:planehom}, i.e.~$|a|=2$, $|b|=1$. The same proof carries over to give $\partial a=\partial b=0$.
\end{proof}

\begin{rem}
The above lemma shows that the subspace of $HLC_\bullet(L_{knot} \cup F;\Z_2)$ spanned by the mixed chords is isomorphic to $H_{\bullet+1}(S^1;\Z_2)$. This graded vector space is, in turn, is isomorphic to the Morse homology of a generic function in the generating family considered above (with shifted degrees).
\end{rem}

\begin{cor}
For every isotopy class of $L_{knot}$, its link with $F$ has a mixed Reeb chord.
\end{cor}

\begin{proof}
If $L_{knot}$ and $F$ could be unlinked, i.e.~if $L_{knot}$ could be isotoped so that the link carries no mixed Reeb chords, then we would have
\[HLC_\bullet(L_{knot} \cup F;\Z_2) \simeq HLC_\bullet(L_{knot};\Z_2) \oplus HLC_\bullet(F;\Z_2) \simeq \Z_2c,\]
where we again have chosen the unique (trivial) augmentations. This leads to a contradiction.
\end{proof}

\begin{proof}[Proof of Theorem \ref{thm:knottedsphere}]
Suppose that $L_{knot}$ is Legendrian isotopic to $L_{std}$. It would then be possible to unlink $L_{knot}$ and $F$, contradicting the previous corollary.

We now show that $\Pi_L(L_{knot}) \subset D^*S^2 \subset S^2 \times S^2$, where $D^*S^2$ denotes the unit disk bundle, cannot be Hamiltonian isotoped to $\Pi_L(L_{std})$ inside $S^2 \times S^2$. We use the fact that the torus considered in \cite{NonDisp}, and fibre-wise rescalings of it, are non-displaceable in $S^2 \times S^2$, as is shown in \cite{ToricDegeneration}. After Hamiltonian isotopy, such a torus can be placed in an arbitrarily small neighbourhood of $\Pi_L(L_{knot})$.

A Hamiltonian isotopy of $S^2 \times S^2$ mapping $\Pi_L(L_{knot})$ to a small Darboux chart, would do the same with a non-displaceable torus sufficiently close to it, which leads to a contradiction.
\end{proof}

\bibliographystyle{gtart}
\bibliography{references}

\end{document}